\documentclass[12pt]{article}
\usepackage{amsmath,amsfonts,amsthm,amssymb}
\usepackage[square]{natbib}
\usepackage[utf8]{inputenc}
\usepackage[numbib]{tocbibind}
\usepackage{microtype,ulem}
\usepackage{hyperref}
\usepackage{xcolor,caption,subcaption}
\hypersetup{
    colorlinks,
    linkcolor={red!50!black},
    citecolor={blue!50!black},
    urlcolor={blue!80!black}
}
\usepackage[inner=2.6cm,outer=1.8cm]{geometry}
\usepackage{tikz,pgfplots}
\usetikzlibrary{external,arrows.meta}
\newcommand{\R}{\mathbb{R}}
\newcommand{\E}{\mathbb{E}}
\newcommand{\C}{\mathbb{C}}
\newcommand{\N}{\mathbb{N}}
\newcommand{\scal}[2]{\left\langle #1,#2 \right\rangle}
\newcommand{\Tr}{\mbox{\rm Tr}}

\newtheorem{thm}{Theorem}[section]
\newtheorem{lem}[thm]{Lemma}
\newtheorem{prop}[thm]{Proposition}

\newtheorem{rem}[thm]{Remark}
\newtheorem{defi}[thm]{Definition}
\newtheorem{hyp}{Assumption}

\title{Lecture notes on non-convex algorithms for low-rank matrix recovery}
\date{}
\author{Irène Waldspurger\thanks{CNRS, Université Paris Dauphine, Inria Mokaplan, France
  (\texttt{waldspurger@ceremade.dauphine.fr})}}

\begin{document}

\maketitle

\begin{abstract}
  Low-rank matrix recovery problems are inverse problems which naturally arise in various fields like signal processing, imaging and machine learning. They are non-convex and NP-hard in full generality. It is therefore a delicate problem to design efficient recovery algorithms and to provide rigorous theoretical insights on the behavior of these algorithms. The goal of these notes is to review recent progress in this direction for the class of so-called ``non-convex algorithms'', with a particular focus on the proof techniques.

  Although they aim at presenting very recent research works, these notes have been written with the intent to be, as much as possible, accessible to non-specialists.
\end{abstract}

These notes were written for an eight-hour lecture at Collège de France. The original version, in French, is available online\footnote{\url{https://www.ceremade.dauphine.fr/~waldspurger/}} and the videos of the lecture can be found on the Collège de France website\footnote{\url{https://www.college-de-france.fr/site/cours-peccot/p10541769392068432_content.htm}}.

The beginning takes inspiration from the review articles \citep*{davenport_romberg} and \citep*{chen_chi}.

\section{Introduction}

In these notes, we consider a family of problems, called \textit{low-rank matrix recovery problems}, which have various applications in data analysis or imaging sciences. The goal is to give an overview of recent theoretical results on a family of algorithms which can be used to solve them, namely \textit{non-convex algorithms}.

The first part of the introduction (Subsection \ref{ss:def_and_ex}) defines low-rank matrix recovery problems and presents examples. The second one (Subsection \ref{ss:goal}) explains what non-convex algorithms are, which theoretical properties of these algorithms will be of interest for us, and why it is important to understand them. The third one (Subsection \ref{ss:outline}) presents the organization of these notes.

\subsection{Low-rank matrix recovery: definition and examples\label{ss:def_and_ex}}

Let $\mathbb{K}$ be either $\R$ or $\C$, and let $n_1,n_2$ be positive integers.

A \textit{low-rank matrix recovery problem} informally consists in recovering an unknown matrix $X^s\in\mathbb{K}^{n_1\times n_2}$ from
\begin{itemize}
\item some ``simple'' information about $X^s$, modeled by the property ``$X^s\in\mathcal{E}$'' for a known subset $\mathcal{E}$ of $\mathbb{K}^{n_1\times n_2}$; this information is often a set of linear measurements over the coefficients of $X^s$, so that $\mathcal{E}$ is an affine subspace of $\mathbb{K}^{n_1\times n_2}$;
\item the knowledge that $X^s$ has low rank, that is, $\mathrm{rank}(X^s)\ll \min(n_1,n_2)$; sometimes, the rank is exactly known, sometimes not.
\end{itemize}
It is often solved by looking for the element in $\mathcal{E}$ with minimal rank:
\begin{equation}\label{eq:min_rank}
  \mbox{minimize }\mathrm{rank}(X)\mbox{ for }X\in\mathcal{E}.
  \tag{min-rank}
\end{equation}

In the following paragraphs, we describe three important examples of low-rank matrix recovery problems.

\subsubsection{First example: matrix completion}

In matrix completion problems, the ``simple'' information is the knowledge of some coefficients of $X^s$, that is $X^s_{i,j}$ for all pairs $(i,j)$ belonging to some set $\Omega\subset\{1,\dots,n_1\}\times\{1,\dots,n_2\}$. This leads to the following minimization problem:
\begin{align*}
  \mbox{minimize }&\mathrm{rank}(X), \\
  \mbox{with }&X_{i,j}=X^s_{i,j},\quad \forall (i,j)\in\Omega.
\end{align*}

One of the reasons which made matrix completion popular is that it serves as a modelization of the \textit{Netflix} problem: for any $k,l$, the coefficient $X^s_{k,l}$ represents the rating user $k$ would give to movie $l$ if s/he watched this movie. The known coefficients $X^s_{i,j},(i,j)\in\Omega$ are the ratings given by users to movies they indeed watched. The goal is to determine the ``not yet given'' ratings. The low-rank assumption models the similarities between movies, as well as between users.

\subsubsection{Second example: phase retrieval}

At first sight, phase retrieval is not a matrix recovery problem: the unknown object is not a matrix, but a vector $x^s\in\mathbb{C}^n$, for some $n\in\N^*$. For some \textit{measurement vectors} $v_1,\dots,v_m\in\C^n$, we have access to
\begin{equation*}
  |\scal{x^s}{v_k}|,\quad
  \mbox{for all }k\leq m.
\end{equation*}
Here, ``$|.|$'' denotes the standard complex modulus, and ``$\scal{.}{.}$'' the usual Hermitian product. The goal is to recover $x^s$. Since, for any $\phi\in\mathbb{R}$,
\begin{equation*}
  |\scal{e^{i\phi}x^s}{v_k}| = |e^{i\phi}|\,|\scal{x^s}{v_k}|=|\scal{x^s}{v_k}|,\quad
  \forall k=1,\dots,m,
\end{equation*}
exact identification of $x^s$ is actually never possible. We therefore only try to recover $x^s$ \textit{up to a global phase}.

This problem is called ``phase retrieval'' because finding the phases of $\scal{x^s}{v_1},\dots,\scal{x^s}{v_m}$ suffices to solve it (if the phases are known, recovering $x^s$ simply amounts to solving a linear system).

Phase retrieval problems naturally appear in optics, and have for this reason been studied since the 1950s. This is in part due to the fact that electromagnetic waves can be modeled as complex-valude functions, whose modulus is much easier to measure than the phase. Detailed explanations can be found in the review article \citep*{schechtman}.

For some families of measurement vectors $v_1,\dots,v_m$, $x^s$ is not uniquely determined by the modulus $|\scal{x^s}{v_1}|,\dots,|\scal{x^s}{v_m}|$, even up to a global phase: another vector $x'$, different from $x^s$ (and from $e^{i\phi}x^s$ for all $\phi\in\R$) can exist such that
\begin{equation*}
|\scal{x^s}{v_k}|=|\scal{x'}{v_k}|,\quad\forall k=1,\dots,m.
\end{equation*}
In this case, reconstructing $x^s$ with certainty is impossible. In the following, we assume that the phase retrieval problems we consider do not suffer from this uniqueness issue: the modulus uniquely determine $x^s$ up to a global phase. This property is notably known to hold for ``generic''\footnote{We say that a property holds true for \textit{generic} vectors if it is satisfied by all $(v_1,\dots,v_m)$ in $(\C^n)^m-\Gamma$, for some subset $\Gamma$ of $(\C^n)^m$ with zero Lebesgue measure.} families of measurement vectors when $m\geq 4n-4$ \citep*{balan_painless}.

A phase retrieval problem can be turned into an equivalent low-rank matrix recovery problem with the change of variable $X^s=x^s(x^s)^* = (x^s_i\overline{x^s_j})_{1\leq i,j\leq n}\in\C^{n\times n}$, called \textit{lifting} \citep*{chai,candes2}. Indeed, for all $x\in\C^n,k\in\{1,\dots,m\}$,
\begin{align*}
  \left(|\scal{x}{v_k}| = |\scal{x^s}{v_k}|\right)
  & \iff \left(|\scal{x}{v_k}|^2 = |\scal{x^s}{v_k}|^2\right) \\
  & \iff \left(\sum_{1\leq i,j\leq n} \overline{x_i}x_j v_{k,i}\overline{v_{k,j}}
    =|\scal{x^s}{v_k}|^2 \right) \\
  & \iff \left( \scal{xx^*}{v_kv_k^*} = |\scal{x^s}{v_k}|^2 \right).
\end{align*}
(In the last line, the notation ``$\scal{.}{.}$'' refers to the usual scalar product over $\mathcal{H}_n(\C)$, the set of Hermitian $n\times n$ matrices: $\scal{A}{B}=\Tr(A^* B)$, with $A^*$ the conjugate transpose of $A$.)

Because a matrix $X\in\mathcal{H}_n(\C)$ can be written as $X=xx^*$ for some $x\in\C^n$ if and only if
\begin{equation*}
  X\succeq 0\quad\mbox{and}\quad\mathrm{rank}(X)=1,
\end{equation*}
the following two problems are equivalent:
\begin{align*}
  \mbox{find }&x\in\C^n \\
  \mbox{such that }&|\scal{x}{v_k}|=|\scal{x^s}{v_k}|,\quad\forall k\leq m;
\end{align*}
\begin{align*}
  \mbox{find }&X\in\mathcal{H}_n(\C) \\
  \mbox{such that }&\scal{X}{v_kv_k^*}=|\scal{x^s}{v_k}|^2,\quad\forall k\leq m, \\
              &X\succeq 0, \\
              & \mathrm{rank}(X)=1.
\end{align*}
The second one is a low-rank matrix recovery problem.

\subsubsection{Third example: phase synchronization\label{sss:phase_sync}}

A phase synchronization problem consists in recovering $n$ complex numbers with unit modulus, $z^s_1,\dots,z^s_n$, from the approximate knowledge of $z^s_k\overline{z^s_l}$, for all $k,l\in\{1,\dots,n\}$, that is from
\begin{equation*}
  C_{k,l} = z^s_k\overline{z^s_l} + w_{k,l},\quad
  \forall k,l\in\{1,\dots,n\},
\end{equation*}
where $w_{k,l}$ is a random noise.

As in phase retrieval, it is never possible to reconstruct the $z^s_k$ more precisely than up to a global phase.

Let us note that reconstruction from the \textit{exact} knowledge of the $z^s_k\overline{z^s_l}$ would be easy: because of the global phase ambiguity, we could assume $z^s_1=1$ and then, for any $k=2,\dots,n$,
\begin{equation*}
z^s_k\overline{z^s_1} = z^s_k,
\end{equation*}
so we have access to $z_k^s$.
What makes the problem difficult is the noise.

Phase synchronization can be seen as a low-rank matrix recovery problem, for the same reason as phase retrieval. The change of variable $Z^s=z^s(z^s)^*$ indeed allows to rewrite it under the equivalent form
\begin{align*}
  \mbox{find }& Z\in\mathcal{H}_n(\C) \\
  \mbox{such that }&Z_{k,l}\approx C_{k,l},\quad \forall k,l\leq n, \\
              & Z\succeq 0, \\
              & \mathrm{rank}(Z)=1.
\end{align*}

This problem is motivated by applications like server synchronization in computer networks. It also serves as a simplified model for the rotation synchronization problem, where one must identify $n$ rotations of $\R^3$, $R_1,\dots,R_n$, from crude estimations of the $R_kR_l^{-1}$. Rotation synchronization is notably important for cryo-electron microscopy; precise references can be found in \citep*{bandeira_boumal_singer}.

\subsection{What are these notes about?\label{ss:goal}}

When confronted with a specific low-rank matrix recovery problem (as happens for most inverse problem), one has to face the following three questions:
\begin{itemize}
\item Identifiability: does the reformulation of the problem under form \eqref{eq:min_rank} really allow to recover our matrix of interest $X^s$? More formally, does Problem \eqref{eq:min_rank} have a unique solution, which is precisely $X^s$?
\item Stability: if there are a few errors in the available information (that is, we only have approximate knowledge of the set $\mathcal{E}$), does Problem \eqref{eq:min_rank} still allow to recover $X^s$, or at least a matrix close to $X^s$?
\item Algorithms: which algorithms are able to solve Problem \eqref{eq:min_rank} as precisely and fast as possible?
\end{itemize}
These three questions are interesting and difficult. In these notes, we only focus on the third one. And since it is only useful to numerically solve Problem \eqref{eq:min_rank} when it allows to identify $X^s$, we will implicitely assume that we are only facing low-rank recovery problems for which identifiability holds.

More specifically, we will be interested in algorithms for which rigorous correctness guarantees are available. An ideal \textit{correctness guarantee} is a statement of the form
\begin{center}
``For all instances in some class $A$ of low-rank recovery problems, algorithm $B$ outputs the correct matrix $X^s$.''
\end{center}
However, most natural classes of low-rank matrix recovery problems are NP-hard in full generality \citep*{hardt14,fickus}. Consequently, if we restrict ourselves to algorithms running in a reasonable amount of time, statements of the above form are most of the time too strong to be true. Therefore, we will mostly discuss weaker forms of correctness guarantees:
\begin{center}
``If we run it on a random element of some class $A$ of low-rank recovery problems, algorithm $B$ outputs the correct matrix $X^s$ with probability close to $1$.''
\end{center}
This is the subject of these notes: the goal is to describe algorithms for which a guarantee of this form can be proved, and to explain the associated proof techniques.

Actually, we will restrict ourselves to \textit{non-convex} algorithms with rigorous correctness guarantees. What are non-convex algorithms? Very broadly speaking, the last decades of research in optimization have allowed for the development of efficient algorithms able to solve with arbitrary precision problems of the form ``minimize $f(x)$ over all $x\in\mathcal{E}$'' when $\mathcal{E}$ is a convex set and $f:\mathcal{E}\to\R$ a convex function. This does not apply to Problem \eqref{eq:min_rank} because the rank is not a convex function. Designing an algorithm applicable to Problem \eqref{eq:min_rank} requires overcoming the non-convexity issue. Two strategies exist.
\begin{itemize}
\item Convex methods: their principle is to approximate the non-convex problem with a convex one, different but with the same minimizer, and solve the convex one. These methods generally work well (in the sense that they correctly recover the solution $X^s$). Moreover, powerful analysis techniques exist to establish correctness guarantees. Since their introduction around 2010, they have therefore been the subject of intense research. On the negative side, they are often impractical because of their large computational cost.
\item Non-convex methods: these algorithms simply ignore the non-convexity of Problem \eqref{eq:min_rank}, and try to solve it with relatively simple strategies, often coming from the field of convex optimization. Most of them are iterative: starting from an initial estimate of the solution, they progressively refine it with natural heuristics. Much older and more widely-used than convex methods, these algorithms have the advantage to be fast. They can a priori fail if they reach an estimate of the solution which heuristics are not able to refine further although it is not the solution; this is called ``reaching a \textit{critical point}''. However, in practice, they often work well.
\end{itemize}
The two strategies have led to the development of useful algorithms and to interesting theoretical analyses. But since it is not possible to cover both in a single eight-hour lecture, the notes focus on the second one only.

Let us stress that the development of rigorous theoretical guarantees for non-convex methods has been the subject of many publications in the last five years. These notes are not an attempt to exhaustively describe all of them. As the main aim is to explain the most successful proof techniques of the domain, we will pick one or two representative articles for each family of techniques, and describe only these, but with some amount of detail.

\subsection{Outline\label{ss:outline}}

In Section \ref{s:convex}, we give a brief overview of convex methods. These methods have played a crucial role in the development of low rank recovery algorithms with correctness guarantees. As far as I know, they are the first ones for which convincing guarantees were established in relatively general settings, and they have served as examples for the theoretical study of other algorithms. It is therefore good to have an idea of what they are and how they work even if they are not our main subject of interest. In addition, convex and non-convex algorithms are not as disjoint categories as they seem, and it has recently been realized that establishing correctness guarantees for a convex method could help establishing guarantees for a non-convex one or vice-versa \citep*{zhong_boumal,chen_chi_fan_ma_yan}.

The short Section \ref{s:non_convex} defines non-convex methods, and provides an example of such a method in the context of phase retrieval.

In Section \ref{s:no_bad_critical}, we describe a first set of techniques which can prove the correctness of some non-convex low-rank recovery algorithms. As said a few paragraphs ago, the main reason why non-convex algorithms can fail is the presence of \textit{critical points}. But, in some settings, it turns out that critical points do simply not exist, which almost automatically implies that non-convex methods succeed. We provide examples of settings where this property holds, and explain how it can be proved.

In Section \ref{s:leave_one_out}, we consider the situation where non-convex methods succeed but critical points exist. This situation is much more delicate than the previous one: establishing correctness guarantees then requires to carefully study the trajectory of iterates of the algorithm and to show that it does not come close to one of the critical points. Few technical tools exist for this. The main one is \textit{leave-one-out}; we present it through the example of the \textit{generalized power method} for phase synchronization \citep*{zhong_boumal}.

Finally, Section \ref{s:burer_monteiro} presents a different class of results. While Sections \ref{s:no_bad_critical} and \ref{s:leave_one_out} explain how to prove the correctness of specific algorithms, applied to specific low-rank recovery problems, Section \ref{s:burer_monteiro} presents a family of algorithms which can be applied to a large class of low-rank recovery problems, and explains why some algorithms in this family can be proved to (almost) always succeed. This provides much more general correctness guarantees than in the two preceding sections. The price to pay is that the algorithms for which these guarantees hold are not as satisfying, from a computational point of view, as the algorithms of Sections \ref{s:no_bad_critical} and \ref{s:leave_one_out}.

\subsection{Notation}

For any vector $x\in\R^n$ or $\C^n$, we denote $||x||$ its $\ell^2$-norm.

Let $X\in\mathbb{K}^{n_1\times n_2}$ be a matrix. We denote $\lambda_1(X),\dots,\lambda_{\min(n_1,n_2)}(X)$ its singular values, that is the nonnegative real numbers such that there exist orthogonal matrices $U\in\mathbb{K}^{n_1\times n_1},V\in\mathbb{K}^{n_2\times n_2}$ for which
\begin{equation*}
  X = U \begin{pmatrix}\lambda_1(X)&&\\&\ddots&\\&&\lambda_{\min(n_1,n_2)}(X) \\
    0&\dots&0\\ \vdots&&\vdots \\ 0&\dots&0
  \end{pmatrix}V.
\end{equation*}
The nuclear norm of $X$ is
\begin{equation*}
||X||_* = \sum_{s=1}^{\min(n_1,n_2)}\lambda_s(X)
\end{equation*}
and its operator norm is
\begin{equation*}
|||X||| = \max_{s=1,\dots,\min(n_1,n_2)}\lambda_s(X).
\end{equation*}

When $\mathbb{K}=\C$, we denote $X^*$ the conjugate transpose of $X$.

We denote $\mathcal{S}_n(\R)$ the set of real symmetric $n\times n$ matrices, $\mathcal{S}_n^+(\R)$ the set of positive semidefinite $n\times n$ matrices, and $\mathcal{H}_n(\C)$ the set of Hermitian $n\times n$ complex matrices.

\section{Convex methods\label{s:convex}}

\subsection{Principle\label{ss:principle}}

The principle of convex methods is to approximate the non-convex problem \eqref{eq:min_rank} with a convex one. The approximation typically consists in replacing the rank (which is a non-convex function) with the nuclear norm (which is a convex one). This replacement is motivated by the fact that $||.||_*$ is the convex envelope of the rank on $\{X\in\mathbb{K}^{n_1\times n_2},|||X|||\leq 1\}$, that is, for all $X\in\mathbb{K}^{n_1\times n_2}$ such that $|||X|||\leq 1$,
\begin{equation*}
||X||_* = \max\{F(X),F:\{X\in\mathbb{K}^{n_1\times n_2},|||X|||\leq 1\} \to \R\mbox{ is convex}, F\leq\mathrm{rank}\}.
\end{equation*}

The problem \eqref{eq:min_rank}
\begin{equation*}
\mbox{minimize }\mathrm{rank}(X)\mbox{ for }X\in\mathcal{E}
\end{equation*}
is thus typically replaced with
\begin{equation}\label{eq:convex_relaxation}
  \mbox{minimize }||X||_* \mbox{ for }X\in\mathcal{E}.
  \tag{Convex relaxation}
\end{equation}
This latter problem is convex if $\mathcal{E}$ is a convex set. It is therefore easier to numerically solve. In addition, although it is a priori only an approximation of Problem \eqref{eq:min_rank}, it often happens to have the same solution, so that solving the convex problem also solves the non-convex one.

This idea was introduced in \citep*{recht_fazel_parrilo}. It echoes the theory of \textit{compressed sensing} for the reconstruction of sparse signals.

\subsection{Convex methods for matrix completion\label{ss:convex_completion}}

Conditions under which the solution of Problem \eqref{eq:convex_relaxation} is the same as the one of Problem \eqref{eq:min_rank} have been established for all low-rank matrix recovery problems described in the introduction. In the case of matrix completion, for instance, Problem \eqref{eq:convex_relaxation} is instanciated as
\begin{align}
  \mbox{minimize }&||X||_* \nonumber\\
  \mbox{with }&X_{i,j} = X^s_{ij},\quad \forall (i,j)\in\Omega.
                \label{eq:convex_completion}
\end{align}
It has notably been analyzed in \citep*{candes_recht, candes_tao, gross}. The simplest situation in which one can guarantee that its solution is indeed $X^s$ (the solution of Problem \eqref{eq:min_rank}) is described by the following theorem.

\begin{thm}[\citet*{chen15}]
Let $r$ be the rank of $X^s\in\R^{n_1\times n_2}$ and $\mu_0$ its incoherence\footnote{Let $U\in\R^{n_1\times r}$ (respectively $V\in\R^{n_2\times r}$) be a matrix whose columns form an orthonormal basis of $\mathrm{Range}(X^s)$ (respectively $\mathrm{Range}(X^{sT})$). Let $U_1,\dots,U_{n_1},V_1,\dots,V_{n_2}$ be the lines of $U,V$. The incoherence $\mu_0$ is defined as $\mu_0=\max\left(
      \frac{n_1}{r}||U_1||^2,\dots,\frac{n_1}{r}||U_{n_1}||^2,
      \frac{n_2}{r}||V_1||^2,\dots,\frac{n_2}{r}||V_{n_2}||^2
    \right)$. Intuitively, it quantifies the maximal alignement between a vector of $\mathrm{Range}(X^s)$ (or $\mathrm{Range}(X^{sT})$) and the canonical basis.}.

  Let us assume that $\Omega$ is chosen at random, and contains each of the $n_1n_2$ pairs $(i,j),i\leq n_1,j\leq n_2$ with probability $p>0$, independently one from each other.

  There exist constants $C,c>0$ such that, if
  \begin{equation*}
    n_1n_2p \geq C \mu_0 r (n_1+n_2)\log^2(\max(n_1,n_2)),
  \end{equation*}
  then, with probability at least $1-\frac{1}{\max(n_1,n_2)^c}$, the solution of Problem \eqref{eq:convex_completion} is $X^s$.
\end{thm}

In other words, as soon as the number of revealed coefficients of $X^s$ is of order $r(n_1+n_2)$ (neglecting $\mu_0$ and the logarithmic term), it is possible to recover $X^s$ by solving Problem \eqref{eq:convex_completion}, with high probability. The number of degrees of freedom of a rank $r$ matrix is of order $r(n_1+n_2)$. The convex algorithm therefore succeeds with (almost) the smallest possible amount of information.

\subsection{Convex methods for phase retrieval}

We recall the ``low-rank matrix recovery'' formulation of phase retrieval problems:
\begin{align*}
  \mbox{find }&X\in\mathcal{H}_n(\C) \\
  \mbox{such that }&\scal{X}{v_kv_k^*}=|\scal{x^s}{v_k}|^2,\quad\forall k\leq m, \\
              &X\succeq 0, \\
              & \mathrm{rank}(X)=1.
\end{align*}
Here, as previously, $x^s$ denotes the vector we want to recover, i.e. the true solution of the phase retrieval problem. The above matricial formulation has been obtained from the initial vectorial one through the change of variable $X=xx^*$, hence its solution must be $X^s=x^sx^{s*}$.

To approximate this non-convex problem with a convex one, we follow the method described in Subsection \ref{ss:principle}: we replace the constraint ``$\mathrm{rank}(X)=1$'' with the constraint ``$||X||_*$ is minimal''. For any $X\succeq 0$, it holds $||X||_*=\Tr(X)$, which leads to the following convex problem (called \eqref{eq:PhaseLift} in \citep*{candes2}):
\begin{align}
  \mbox{minimize }&\Tr(X)\nonumber \\
  \mbox{with }&\scal{X}{v_kv_k^*}=|\scal{x^s}{v_k}|^2,\quad\forall k\leq m,
                \tag{PhaseLift}\label{eq:PhaseLift} \\
  & X\succeq 0.\nonumber
\end{align}
Several works have established correctness guarantees for \eqref{eq:PhaseLift} (that is, have proved, in various situations, that the solution of \eqref{eq:PhaseLift} was $X^s=x^sx^{s*})$. Most of them apply in the setting where $v_1,\dots,v_m$ are chosen according to independent normal laws:
\begin{equation*}
v_1,\dots,v_m\overset{iid}{\sim} \mathcal{N}(0,I_n).
\end{equation*}
This assumption is somewhat unrealistic: in most applications, the measurement vectors are not random; even when they are, they seldom follow a normal law. However, this assumption turns out to simplify the computations a lot, and it is difficult to get rid of it.

Under this assumption, the strongest available guarantees have been proved in \citep*{candes_li}, and are summarized in the following theorem.
\begin{thm}\label{thm:candes_li}
  There exist constants $C,c>0$ such that, for any $n,m\in\N^*$, if $m\geq Cn$, it holds with probability at least $1-e^{-cm}$ that, whatever $x^s\in\C^n$, the solution $\hat X$ of Problem \eqref{eq:PhaseLift} is
  \begin{equation*}
    \hat X = x^s(x^s)^*.
  \end{equation*}
\end{thm}
Let us note that, when $m<n$, $x^s$ is never uniquely determined by $|\scal{x^s}{v_1}|,\dots,|\scal{x^s}{v_m}|$ up to a global phase, hence the phase retrieval problem is intrinsically not solvable. The theorem therefore means that, with high probabiliy, \eqref{eq:PhaseLift} solves phase retrieval problems with an optimal (up to a multiplicative constant) number of measurements.

\begin{proof}[Sketch of proof, adapted from \citep*{chen_chi_goldsmith}] We define
  \begin{gather*}
    \mathcal{A}:X \in \mathcal{H}_n(\C) \to (\scal{X}{v_kv_k^*})_{k=1,\dots,m}\in\R^m;\\
    \tilde{\mathcal{A}}:X\in\mathcal{H}_n(\C)
    \to (\mathcal{A}(X)_1-\mathcal{A}(X)_2,\mathcal{A}(X)_3-\mathcal{A}(X)_4,\dots)
    \in\R^{\lfloor m/2\rfloor}.
  \end{gather*}
  Here, for any $k$, $\mathcal{A}(X)_k$ denotes the $k$-th coordinate of $\mathcal{A}(X)$. Defining $\tilde{\mathcal{A}}$ is necessary to make the proof correct but, in order to grasp the main ideas of the sketch, one can simply imagine that $\tilde{\mathcal{A}}=\mathcal{A}$.

  \begin{prop}
    There exist $c,C,d,D>0$ such that, for any $r\in\{1,\dots,n\}$, the following property is true: if $m\geq Cnr$, then, with probability at leat $1-e^{-cm}$,
    \begin{equation}\label{eq:RIP}
      \tag{RIP}
      d ||X||_F \leq \frac{1}{m}||\tilde{\mathcal{A}}(X)||_{\ell^1}\leq D ||X||_F.
    \end{equation}
    for all matrices $X\in\mathcal{H}_n(\C)$ with rank at most $r$.

    (Here, $||X||_F$ is the Frobenius norm of $X$: $||X||_F = \left(\sum_{1\leq i,j\leq n}|X_{i,j}|^2\right)^{1/2}$.)
  \end{prop}
  When \eqref{eq:RIP} holds, we say that $\tilde{\mathcal{A}}$ satisfies a \textit{Restricted Isometry Property} in $\ell^1$ / Frobenius norms (notion which echoes the theory of compressed sensing again). When it holds, this property implies that $X=x^s(x^s)^*$ is the unique solution of Problem \eqref{eq:PhaseLift}. Indeed, if we denote $\hat X$ the solution of the problem, we must have
  \begin{gather*}
    \Tr(\hat X)\leq \Tr(x^s(x^s)^*), \\
    \scal{\hat X}{v_kv_k^*}=|\scal{x^s}{v_k}|^2
    =\scal{x^s(x^s)^*}{v_kv_k^*},\quad \forall k=1,\dots,m, \\
    \hat X\succeq 0,
  \end{gather*}
  which implies, for $H\overset{def}{=}\hat{X} - x^s(x^s)^*$:
  \begin{subequations}
    \begin{gather}
      \Tr(H)\leq 0, \label{eq:propH1} \\
      \tilde{\mathcal{A}}(H)=0, \label{eq:propH2} \\
      x^s(x^s)^* + H \succeq 0. \label{eq:propH3}
    \end{gather}
  \end{subequations}
  The following proposition concludes.
  \begin{prop}
    If Property \eqref{eq:RIP} holds for some $r > \frac{4D^2}{d^2}+2$, then, whatever $x^s\in\C^n$, there does not exist $H\in\mathcal{H}_n(\C)$ a non-zero matrix for which Properties \eqref{eq:propH1}, \eqref{eq:propH2} and \eqref{eq:propH3} are true.
  \end{prop}
  A very vague expanation of this proposition is that, if $H$ is such that $\Tr(H)\leq 0$ and $x^s(x^s)^*+H\succeq 0$, then $H$ is ``somewhat close to $-\lambda x^s(x^s)^*$'' for some $\lambda >0$. As a consequence, even if the rank of $H$ may be larger than $r$, Property \eqref{eq:RIP} allows to show that
  \begin{equation*}
    \frac{1}{m}||\tilde{\mathcal{A}}(H)||_{\ell^1} \geq d' ||H||_F
  \end{equation*}
  for some constant $d'>0$. The equality $\tilde{\mathcal{A}}(H)=0$ therefore implies that $H=0$.
\end{proof}

Many extensions of this result exist: it is possible to analyze the stability of \eqref{eq:PhaseLift} to noise, to consider other random distributions for the measurement vectors than a normal law, wonder how to exploit additional knowledge we might have on $x^s$ ... A non-exhaustive list of references is \citep*{candes_li2, gross, Li12}.

Other convex phase retrieval methods can also be designed. In particular, it turns out that the requirement, in \eqref{eq:PhaseLift}, that the trace of $X$ is minimal, is not stricly necessary. Surprisingly, the following problem
\begin{align}
  \mbox{find }&X\in\mathcal{H}_n(\C) \nonumber \\
  \mbox{such that }&\scal{X}{v_kv_k^*}=|\scal{x^s}{v_k}|^2,\quad \forall k\leq m,
                     \tag{Weak-PhaseLift}\label{eq:WeakPhaseLift} \\
              &X\succeq 0\nonumber
\end{align}
satisfies almost the same correctness guarantees as the ones stated in Theorem \ref{thm:candes_li} \citep*{demanet}. With a suitable change of variable, Problem \eqref{eq:WeakPhaseLift} can be transformed into another problem, called \eqref{eq:PhaseCut},
\begin{align}
  \mbox{minimize }&\Tr(MU) \nonumber\\
  \mbox{with }&U\in\mathcal{H}_m(\C), \nonumber \\
                  &U_{k,k} = 1,\quad \forall k\leq m,\tag{PhaseCut}
                    \label{eq:PhaseCut} \\
  & U\succeq 0,\nonumber
\end{align}
where $M$ is an explicit matrix which depends on the $v_k$ and $|\scal{x^s}{v_k}|$. The \eqref{eq:PhaseCut} formulation has some advantages over the \eqref{eq:PhaseLift} one, from an algorithmic point of view. It also seems to be a bit more robust to noise in some settings \citep*{maxcut}.

\subsection{Limits of convex methods}

To summarize, convex methods allow solving low-rank matrix recovery problems with an essentialy minimal number of measurements. Their algorithmic side is relatively well understood; many numerical solvers have been developed for problems of the form \eqref{eq:convex_relaxation}.

This idyllic picture must however be toned down: convex methods suffer from an intrinsic ``dimensionality issue'', which severely limits the computational efficiency of numerical solvers, and hence the practical applicability of these methods.

A matrix $X^s\in\mathbb{K}^{n_1\times n_2}$, with rank $r$, can be written as
\begin{equation*}
X^s = U V^T
\end{equation*}
for some $U\in\mathbb{K}^{n_1\times r},V\in\mathbb{K}^{n_2\times r}$. When the rank $r$ of the solution is known, the number of ``degrees of freedom'' of a low-rank recovery problem is therefore of order $(n_1+n_2)r$\footnote{slightly smaller actually, since the matrices $U,V$ in the decomposition $X^s=UV^T$ are not unique}. But to solve the convex problem \eqref{eq:convex_relaxation}, one must reconstruct each of the $n_1n_2$ coefficients of $X^s$. Thus, the complexity of numerical solvers tends to be polynomial in $n_1n_2$, which is much larger than what we could have hoped for:
\begin{equation*}
n_1n_2 \gg (n_1+n_2)r.
\end{equation*}

Let us illustrate this phenomenon with two examples of algorithms.
\begin{enumerate}
\item Interior-point methods: these algorithms apply to convex problems of the form
  \begin{align*}
    \mbox{minimize }&\Tr(X) \\
    \mbox{with }&\mathcal{A}(X)=0, \\
    &X\succeq 0,
  \end{align*}
  where $\mathcal{A}:\mathbb{K}^{n\times n}\to\R^m$ is an affine map. When $\mathcal{E}$ is affine, Problem \eqref{eq:convex_relaxation} can be rewritten under this form for $n=n_1+n_2$, with the introduction of properly chosen additional variables.

  Interior-point methods are iterative. The number of iterations necessary to reach a good approximation of the solution is in general moderate but, in full generality, the number of arithmetic operations performed at each iteration \citep*[Page 357]{borchers} is of order
  \begin{equation*}
    (m+n)mn^2.
  \end{equation*}
\item FISTA \citep*{beck_teboulle}: this algorithm approximately solves problems of the form \eqref{eq:convex_relaxation}, also in the case where $\mathcal{E}$ is affine, by replacing them with a ``regularized'' version:
  \begin{equation*}
    \underset{X\in\mathbb{K}^{n_1\times n_2}}{\mbox{minimize}}
    \frac{1}{2}||\mathcal{A}(X)||^2 + \lambda ||X||_*,
  \end{equation*}
  for $\mathcal{A}:\mathbb{K}^{n_1\times n_2}\to\R^m$ an affine map such that $\mathcal{E}=\{X,\mathcal{A}(X)=0\}$, and $\lambda>0$ a parameter.

  This algorithm is also iterative. It requires more iterations than an interior-point method. Each iteration consists in one step of gradient descent on the map $X\to ||\mathcal{A}(X)||^2$ and one application of the so-called \textit{proximal operator} of the nuclear norm. The proximal operator is the most costly part. Heuristics exist to accelerate its computation by exploiting the structure of the problem but, in the most general case, it requires to perform a singular value decomposition, which represents $O(n_1^3)$ operations if $n_1$ and $n_2$ are of the same order.
\end{enumerate}
It is clear that algorithms with these complexities are unapplicable when $n_1,n_2,m$ grow large.

To overcome this dimensionality issue, an option is to improve numerical solvers further, by taking advantage of the fact that convex approximations of low-rank recovery problems are not ``generic'' problems of the form \eqref{eq:convex_relaxation} but (in principle) problems of the form \eqref{eq:convex_relaxation} \textit{with a low-rank solution} \citep*{ding_yurtsever_cevher_tropp_udell, yurtsever_tropp_fercoq_udell_cevher} 

Another option is simply to come back to non-convex methods, either to the traditional heuristics mentioned in Subsection \ref{ss:goal} or to more modern ones. Precisely understanding when non-convex methods succeed and when they fail is still out of reach. But in the last years, a flurry of works has at least allowed to find some situations where non-convex methods succeed and it is moreover possible to rigorously explain it.

\section{Non-convex methods: definition and example\label{s:non_convex}}

\subsection{Informal definition of non-convex methods\label{ss:non_convex_informal}}

A rank $r$ matrix $X\in\mathbb{K}^{n_1\times n_2}$ can always be written under the form
\begin{equation*}
  X = UV^T,\quad\mbox{for some }U\in\mathbb{K}^{n_1\times r},V\in\mathbb{K}^{n_2\times r}.
\end{equation*}
When $X$ is symmetric (respectively Hermitian if $\mathbb{K}=\mathbb{C}$), it can actually be written as $X=UU^T$ (respectively $X=UU^*$) for some $U\in\mathbb{K}^{n_1\times n_1}$. Reconstructing the factors $U$ and $V$ suffices to recover $X$.

\begin{rem}
In the case of low-rank matrix recovery problems coming from phase retrieval or phase synchronization, this property is obvious, since the low-rank matrix $X^s$ (or $Z^s$ in the case of phase synchronization) to recover is precisely defined as $X^s=x^s(x^s)^*$ (respectively $Z^s=z^s(z^s)^*$), with $x^s\in\C^{n\times 1}$ (respectively $z^s\in\C^{n\times 1}$) the solution of the original problem.
\end{rem}

Non-convex algorithms are generally based on the following scheme:
\begin{enumerate}
\item Choice of an initial point $(U_0,V_0)$ (or simply $U_0$ in the symmetric/hermitian case), which can be an estimation of the solution if available or an arbitrary point.
\item Iterative application of some simple heuristics aiming at making $(U_0,V_0)$ closer and closer to a solution (that is, a pair $(U,V)$ such that $X=UV^T$ is the sought-after low-rank matrix $X^s$). This yields a sequence of iterates $(U_t,V_t)_{t\in\N}$ which should ideally converge.
\item Output of $U_TV_T^T$ for some $T\in\N$ large enough.
\end{enumerate}

These methods are widely-used in practice. As we have said, they can fail, but it is empirically observed that they work well in various situations.

\subsection{Example: alternating projections for phase retrieval}

In this subsection, we present an example of non-convex method: the phase retrieval algorithm called \textit{alternating projections}. It is the oldest and most well-known phase retrieval algorithm, sometimes named \textit{error reduction} or \textit{Gerchberg-Saxton} (from the name of the researchers who introduced it \citep*{gerchberg}).

We consider a general phase retrieval problem, in ``vector'' form: we aim at reconstructing $x^s\in\C^n$ from
\begin{equation*}
b_k\overset{def}{=}|\scal{x^s}{v_k}|,\quad \forall k=1,\dots,m.
\end{equation*}
We define $\mathcal{A}:x\in \C^n\to\mathcal{A}(x)=(\scal{x}{v_1},\dots,\scal{x}{v_m})\in\C^m$. If $m\geq n$, $\mathcal{A}$ is injective for generic measurement vectors. Under this assumption, recovering $x^s$ (up to a global phase) is equivalent to recovering $y^s\overset{def}{=}\mathcal{A}(x^s)$ (also up to a global phase). Our phase retrieval problem can thus be reformulated as
\begin{align*}
  \mbox{find }&y\in\C^m, \\
  \mbox{such that }&y\in\mathrm{Range}(\mathcal{A}), \\
  &|y_k|=b_k,\quad \forall k\leq m.
\end{align*}
If we introduce the notation $E=\{y\in\C^m,|y_k|=b_k,\forall k\leq m\}$, this problem can be concisely rewritten as
\begin{align*}
  \mbox{find }&y\in\C^m, \\
  \mbox{such that }&y\in\mathrm{Range}(\mathcal{A})\cap E.
\end{align*}
The alternating projections algorithm is a natural heuristic inspired by this formulation.
\begin{defi}
  Let $T\in\N$ be fixed. The alternating projections method is composed of the following steps.
  \begin{enumerate}
  \item We choose an arbitrary starting point (for our numerical experiments, we will choose $y_0\sim\mathcal{N}(0,y_m)$).
  \item For every $t=1,\dots,T$, we define
    \begin{equation*}
      y_t = P_{\mathrm{Im}(\mathcal{A})}(P_E(y_{t-1})),
    \end{equation*}
    where, for any non-empty closed subset of $\C^m$, $P_S$ denotes the projection onto $S$\footnote{that is to say a (possibly non uniquely defined) function $P_S:\C^m\to S$ such that, for any $y$, $||P_S(y)-y||=\min_{z\in S}||z-y||$}.
  \item We return $y_T$.
  \end{enumerate}
\end{defi}

\begin{rem}
  Although the ``vectorial'' definition we just gave does not make it fully apparent, the alternating projections algorithm can be seen as an instance of the scheme described in Subsection \ref{ss:non_convex_informal}, up to the change of variable
  \begin{equation*}
    U_t = \mathcal{A}^{-1}(y_t),\quad\forall t\geq 1.
  \end{equation*}
  Here, $\mathcal{A}^{-1}$ denotes the inverse of $\mathcal{A}$, when seen as an operator from $\C^n$ to $\mathrm{Range}(\mathcal{A})$.
\end{rem}

\begin{figure}
    \centering
    \captionsetup{width=0.8\textwidth}
\definecolor{mycolor1}{rgb}{0.00000,0.44700,0.74100}%
\definecolor{mycolor2}{rgb}{0.85000,0.32500,0.09800}%
\begin{tikzpicture}

\begin{axis}[%
width=4.602in,
height=3.506in,
at={(0.772in,0.473in)},
scale only axis,
xmin=0,
xmax=8,
ymin=0,
ymax=1.1,
axis background/.style={fill=white},
axis x line*=bottom,
axis y line*=left,
legend style={legend cell align=left, align=left, draw=white!15!black,
  anchor = south west,, at={(0.7,0.05)}},
xlabel={$m/n$},
ylabel={Probabilité},
scale=0.6
]
\addplot [color=mycolor1, line width=2.0pt, mark=x, mark options={solid, mycolor1}]
  table[row sep=crcr]{%
2	0\\
2.5	0.025\\
3	0.199\\
3.5	0.451\\
4	0.658\\
4.5	0.763\\
5	0.856\\
5.5	0.928\\
6	0.95\\
6.5	0.981\\
7	0.976\\
7.5	0.987\\
};
\addlegendentry{Alternating projections}

\addplot [color=mycolor2, line width=2.0pt, mark=+, mark options={solid, mycolor2}]
  table[row sep=crcr]{%
2	0\\
2.5	0\\
3	0\\
3.5	0\\
4	0.041666666666667\\
4.5	0.375\\
5	0.916666666666667\\
5.5	1\\
6	1\\
6.5	1\\
7	1\\
7.5	1\\
};
\addlegendentry{\textit{PhaseLift}}

\end{axis}
\end{tikzpicture}%
    \caption{Success probability for two phase retrieval algorithms, as a function of $m/n$, for $n=40$.\label{fig:proba_rec_phase}}
  \end{figure}
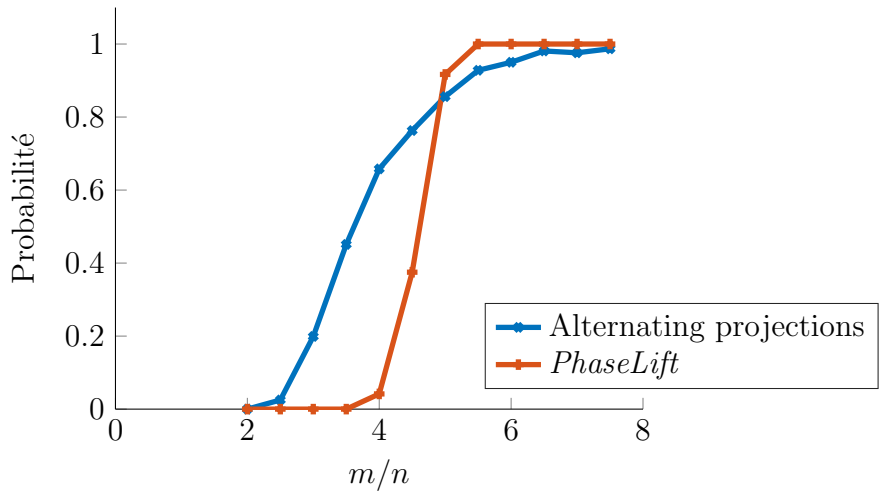

We illustrate the behavior of this method with a numerical experiment in the setting where $v_1,\dots,v_m$ are independently chosen according to a normal distribution:
\begin{equation*}
v_1,\dots,v_m\overset{iid}{\sim}\mathcal{N}(0,I_m).
\end{equation*}
Figure \ref{fig:proba_rec_phase}\footnote{Figures have been generated with the code available at \url{https://www.ceremade.dauphine.fr/~waldspurger/code/non_convex_lecture_notes_figures.zip}, except the red curve on Figure \ref{fig:proba_rec_phase}, generated with \textit{PhasePack} \citep*{phasepack}.} displays the performance curve of the alternating projections method (that is, the probability of exact recovery), as a function of $m/n$, for $n=40$. From this figure, it is tempting to conjecture that, when measurement vectors are normally distributed, alternating projections succeed with high probability as soon as $m\geq Cn$, which would be the same correctness guarantee as the one we have seen for \eqref{eq:PhaseLift} (Theorem \ref{thm:candes_li}). Can we prove it?

We will not answer this question: the alternating projections method is difficult to analyze and the question is still open. However, in the next sections, we will explain how to establish similar guarantees for some other non-convex low-rank recovery algorithms.

\section{When there are no bad critical points\label{s:no_bad_critical}}

The main obstacle to the convergence of non-convex methods is the possible existence of \textit{bad critical points}, at which the heuristic used in the method can get stuck. But, for some heuristics, it is possible to show that such points do not exist, which implies that the algorithms succeed. The goal of this section is to present this proof technique through examples.

Let us first underline that this technique cannot be applied to all non-convex algorithms: most classical heuristics possess bad critical points. It numerically seems, for instance, that the alternating projections method we just discussed possesses bad critical points even when $m/n$ is very large. This does not prevent alternating projections from working well, but makes them out of reach of the proof technique of this section. Actually, among the algorithms to which this technique has been successfully applied, most are not traditional methods used by practioners for decades, but algorithms explicitely designed in order to make them amenable to the proof technique.

The technique can be applied in two slightly different ways:
\begin{itemize}
\item One can
  \begin{enumerate}
  \item first prove that the heuristic used in the algorithm has no bad critical point \textit{in some neighborhood of the solution};
  \item then show that all iterates belong to this neighborhood (which is only possible if the algorithm uses a careful initialization strategy: if the initial point is arbitrary, it has no reason to belong to the neighborhood).
  \end{enumerate}
  This strategy is described in Subsection \ref{ss:no_crit_close}.
\item One can show that the heuristic has no bad critical point at all. This strategy is described in Subsection \ref{ss:no_crit_at_all}.
\end{itemize}
We describe the strategies through examples of phase retrieval algorithms, under the assumption that measurement vectors are normally distributed. However, they apply to many other problems and algorithms. References will be given at the end of Subsections \ref{ss:no_crit_close} and \ref{ss:no_crit_at_all}.

\subsection{No bad critical point close to the solution\label{ss:no_crit_close}}

We illustrate the first method with the study of \textit{Wirtinger Flow}, a phase retrieval algorithm introduced in \citep*{candes_wirtinger}, which consists of the following steps:
\begin{enumerate}
\item Initialization: choice of $x_0$ according to a so-called \textit{spectral method}, which we will describe later.
\item Refinement step: let us define
  \begin{equation*}
    f : x\in\C^n \to \frac{1}{2m}\sum_{k=1}^m(|\scal{x}{v_k}|^2-b_k^2)^2.
  \end{equation*}
  It is a $\mathcal{C}^{\infty}$ non-convex function, whose minima are exactly the solutions of the phase retrieval problem\footnote{For all $x$, $f(x)\geq 0$ and equality is reached if and only if $|\scal{x}{v_k}|=b_k$ for all $k$.}. \\
  For all $t\in\N$, we define $x_{t+1}$ by applying to $x_t$ a gradient descent step over $f$:
  \begin{equation*}
    x_{t+1} = x_t - \mu\nabla f(x_t),
  \end{equation*}
  for some constant $\mu>0$.
\item Output of $x_T$ for some $T$ large enough.
\end{enumerate}

This algorithm obeys the correctness guarantees stated in the following theorem\footnote{In this theorem, the notation ``$\mathrm{dist}$'' is defined as $\mathrm{dist}(x,y)=\min_{\phi\in\R}||x-e^{i\phi}y||$. In the rest of the subsection, we will do as if $\mathrm{dist}(x,y)=||x-y||$, to simplify the explanations, but this is not rigorous.}.

\begin{thm}[{\citep*[Thm 3.3]{candes_wirtinger}}]\label{thm:wirtinger}
  Let the solution $x^s$ of the phase retrieval problem be arbitrary. Let the measurement vectors $v_1,\dots,v_m$ be generated according to independent normal distributions.

  There exist constants $C,c>0$ such that, if
  \begin{equation*}
    C n\log(n) \leq m
  \end{equation*}
  and if $\mu\leq \frac{c}{n}$, then, with probability $1-O\left(\frac{1}{n^2}\right)$,
  \begin{equation*}
    \mathrm{dist}(x^s,x_t) \leq \frac{1}{8}\left(1-\frac{\mu}{4}\right)^{t/2}||x^s||
  \end{equation*}
  for all $t\in\N$. In particular, $x_t\overset{t\to+\infty}{\to} x^s$.
\end{thm}

As said at the beginning of the section, the proof of this theorem proceeds in two steps, whose principles will be described in Paragraphs \ref{sss:first_step} and \ref{sss:second_step}.
\begin{enumerate}
\item We show that the heuristic used in the refinement step has no bad critical point in the ball $B(x^s,||x^s||/8)$. More precisely, we establish that, for any $x\in B(x^s,||x^s||/8)$,
  \begin{equation}\label{eq:contraction_wirt}
    \mathrm{dist}(x^s,x-\mu\nabla f(x)) \leq \rho\,\mathrm{dist}(x^s,x)
  \end{equation}
  for some $\rho = \sqrt{1-\frac{\mu}{4}}\in]0;1[$.
\item We analyze the spectral initialization method and show that
  \begin{equation}\label{eq:init_wirt}
    x_0\in B(x^s,||x^s||/8).
  \end{equation}
\end{enumerate}
These two steps suffice to establish the theorem: starting from Property \eqref{eq:init_wirt} and iteratively applying Equation \eqref{eq:contraction_wirt} immediately proves the main inequality of Theorem \ref{thm:wirtinger}.

\subsubsection{First step\label{sss:first_step}}

We assume to simplify that $||x^s||=1$.

In this paragraph, we explain how to prove Property \eqref{eq:contraction_wirt}, but for some $\rho\in]0;1[$ larger than $\sqrt{1-\frac{\mu}{4}}$. This modification degrades the convergence rate guaranteed by Theorem \ref{thm:wirtinger} but allows for a slightly simpler proof.

It is enough to prove the following properties:
\begin{gather}
  \forall x\in B(x^s,1/8),\quad
  \mathrm{Re}(\scal{x-x^s}{\nabla f(x)}) \geq \alpha \mathrm{dist}(x,x^s)^2,
  \label{eq:contraction1}\\
  \forall x\in B(x^s,1/8),\quad
  ||\nabla f(x)|| \leq \beta \mathrm{dist}(x,x^s),
  \label{eq:contraction2}
\end{gather}
for $\alpha,\beta>0$ some well-chosen values. Indeed, if these inequalities are true, we have for all $x\in B(x^s,1/8)$ that
\begin{align*}
  ||x^s-(x-\mu\nabla f(x))||^2
  & = ||x^s-x||^2 - 2 \mu \mathrm{Re}( \scal{x-x^s}{\nabla f(x)}) + \mu^2 ||\nabla f(x)||^2 \\
  & \leq ||x^s-x||^2 - 2 \mu \alpha ||x-x^s||^2 + \mu^2 \beta^2 ||x-x^s||^2 \\
  & \leq (1-\mu\alpha) ||x^s-x||^2
    \quad\mbox{if }\mu < \frac{\alpha}{\beta^2}.
\end{align*}

The general principle for proving Properties \eqref{eq:contraction1} and \eqref{eq:contraction2} is to compute the explicit expression of $\nabla f$ and use it to write $\mathrm{Re}(\scal{x-x^s}{\nabla f(x)})$ and $||\nabla f(x)||^2$ as a sum of realizations of independent random variables, which can be analyzed with the help of a classical statistical tool: concentration inequalities\footnote{Here is a brief definition of concentration inequalities, for readers who are not familiar with them. In their most basic version, concentration inequalities aim at studying the behavior of a sum of independent random variables
  \begin{equation*}
    Y_1+Y_2+\dots +Y_K.
  \end{equation*}
  Under reasonable assumptions, the sum is close to its expectation with high probability when $K$ is large enough. Concentration inequalities allow to precisely control this closeness, by providing upper bounds for
  \begin{equation*}
    \mathrm{Prob}(Y_1+\dots+Y_K\geq\mathbb{E}(Y_1+\dots+Y_K)+\epsilon)
  \end{equation*}
  for all $\epsilon>0$. The exact form of the upper bound depends on the hypotheses available for the $Y_k$.}.

Let us focus on Property \eqref{eq:contraction1}. Let us consider some $x$ of the form $x=x^s+h$ with $||h||<1/8$. We have
\begin{align*}
  \nabla f(x)&=\frac{1}{m} \sum_{r=1}^m (|\scal{v_r}{x}|^2-b_r^2)v_rv_r^*x \\
  & = \frac{1}{m} \sum_{r=1}^m
    \left(2\mathrm{Re}(\overline{\scal{v_r}{x^s}}\scal{v_r}{h}
    + |\scal{v_r}{h}|^2)\right)\scal{v_r}{x^s+h}v_r,
\end{align*}
which implies that
\begin{align*}
  \mathrm{Re}& (\scal{x-x^s}{\nabla f(x)}) \\
             &= \frac{1}{m}\sum_{r=1}^m
               \left(2 \mathrm{Re}^2(\overline{\scal{v_r}{x^s}}\scal{v_r}{h})
               + 3 \mathrm{Re}(\overline{\scal{v_r}{x^s}}\scal{v_r}{h}) |\scal{v_r}{h}|^2
               + |\scal{v_r}{h}|^4\right) \\
             & \overset{\mbox{\scriptsize def}}{=}
               \frac{1}{m} \sum_{r=1}^m Y_r(h).
\end{align*}
For any fixed $h$, the random variables $Y_1(h),\dots,Y_r(h)$ are independent and identically distributed. One can check that, for all $r$, if we assume (to simplify) that $\scal{x^s}{h}$ belongs to $\R$,
\begin{align*}
  \mathbb{E}(Y_r(h))
  & = 3 \scal{x^s}{h}^2+
    ||h||^2
  + 6 \scal{x^s}{h}||h||^2
  + 2 ||h||^4 \\
  & \geq \frac{||h||^2}{2} \mbox{ if }||h||\leq \frac{1}{8}.
\end{align*}
The above-mentioned concentration inequalities allow to show that
\begin{equation*}
  \frac{1}{m} \sum_{r=1}^m Y_r(h)
  \geq \mathbb{E}\left(\frac{1}{m} \sum_{r=1}^m Y_r(h)\right) - \frac{||h||^2}{4}
\end{equation*}
with probability at least $1-e^{-\gamma m}$ for some constant $\gamma>0$. From this we deduce
\begin{align*}
  \mathrm{Re} (\scal{x-x^s}{\nabla f(x)}) \geq \frac{||h||^2}{4} = \frac{||x-x^s||^2}{4},
\end{align*}
which is the inequality in Property \eqref{eq:contraction1}, with $\alpha=\frac{1}{4}$.

This proof only shows that, for some fixed $x\in B(x^s,1/8)$, the inequality of Property \eqref{eq:contraction1} holds with high probability. It does not show that the inequality holds with high probability for all $x\in B(x^s,1/8)$ at the same time. However, the proof can be extended to all elements $x\in B(x^s,1/8)$ using a very classical probabilistic argument called \textit{$\epsilon$-net}.

\subsubsection{Second step\label{sss:second_step}}

Let us describe the spectral initialization method. As far as I know, this method was first proposed for phase retrieval in \citep*{netrapalli}, before being used in \citep*{candes_wirtinger}. A similar idea could already be found in \citep*{keshavan}, but applied to matrix completion problems.

Let us define the matrix
\begin{equation*}
  M = \frac{1}{m} \sum_{r=1}^m b_r^2 v_r v_r^* = \frac{1}{m} \sum_{r=1}^m |\scal{x^s}{v_r}|^2 v_r v_r^*\in \C^{n\times n}.
\end{equation*}
Informally, in this definition, for each $r$, the rank-$1$ matrix $v_rv_r^*$ appears with a weight proportional to $|\scal{x^s}{v_r}|^2$. Consequently, the more it is aligned with $x^s(x^s)^*$, the more it contributes to $M$, so that $M$ is ``biased'' in the direction of $x^s(x^s)^*$. This reasoning justifies using as an initial point for \textit{Wirtinger Flow}
\begin{equation*}
x_0 = \mbox{main eigenvector of }M.
\end{equation*}

The following lemma establishes precision guarantees for the spectral initialization method.
\begin{lem}
  There exists a constant $C>0$ such that, when $m\geq Cn\log(n)$,
  \begin{equation*}
    \mathrm{dist}(x_0,x^s) \leq \frac{||x^s||}{8}
  \end{equation*}
  with probability $1-O\left(\frac{1}{n^2}\right)$.
\end{lem}

The proof of the lemma relies on the following property, valid with probability $1-O\left(\frac{1}{n^2}\right)$:
\begin{equation*}
|||M-\mathbb{E}(M)|||\leq \delta,
\end{equation*}
where $\delta$ is a constant which can be arbitrarily small if the constant $C$ in the lemma is large enough. This property is proved with the help of concentration inequalities and allows to establish the lemma thanks to the equality
\begin{equation*}
\mathbb{E}(M)=I_n+x^s(x^s)^*.
\end{equation*}

\subsubsection{Related work}

Many other algorithms than \textit{Wirtinger Flow} have been designed, which rely on the same two-step scheme ``spectral initialization + refinement'' and are amenable to a similar theoretical analysis.

Several phase retrieval algorithms notably use this principle, replacing the basic spectral initialization method of \textit{Wirtinger Flow} with a more sophisticated one \citep*{candes_wirtinger2,mondelli} and using another refinement heuristic like ``truncated'' gradient descent \citep*{candes_wirtinger2}, gradient descent on a different cost function than the Wirtinger one (possibly non-smooth, which raises some technical difficulties) \citep*{zhang,wang} or alternating projections \citep*{gerchberg_saxton_short}. These algorithms have correction guarantees similar to the ones described in Theorem \ref{thm:wirtinger}, actually a bit better

For other low-rank matrix recovery problems, we can cite for instance \citep*{jain} for matrix completion and \textit{RIP matrix sensing} (matrix recovery from linear measurements when the measurement operator satisfies a restricted isometry property), \citep*{zhao_wang_liu} for RIP matrix sensing, with a wider range of refinement heuristics than the previous article, \citep*{zheng_lafferty} for matrix completion again and \citep*{chen_wainwright} for more general results, applicable to several matrix recovery problems.

\subsection{No bad critical point at all\label{ss:no_crit_at_all}}

The second proof technique we describe is a variant of the previous one. It consists in showing that the refinement heuristic has no critical point at all (instead of having no critical point \textit{in a neighborhood of the solution}). From a technical point of view, it is in general more involved than the first one: it necessitates an even finer analysis of the functions which come into play. On the positive side, it can apply to conceptually simpler algorithms, closer to practice: with this technique, no sophisticated initialization technique is necessary.

The algorithm we use to illustrate this technique has been proposed in \citep*{sun_qu_wright}. It is identical to \textit{Wirtinger Flow}, except for the initialization:
\begin{enumerate}
\item Initialization: random choice of $x_0$; we may for instance pick $x_0\in B(0,1)$ with uniform probability.
\item Refinement step: gradient descent\footnote{\citet*{sun_qu_wright} actually propose trying to minimize $f$ with another local optimization method than gradient descent, called \textit{Trust-Region}. As the convergence theorem we are going to state holds for both gradient descent and Trust-Region, we use gradient descent.} on the function
  \begin{equation*}
    f : x\in\C^n \quad\to\quad \frac{1}{2m}\sum_{k=1}^m(|\scal{x}{v_k}|^2-b_k^2)^2.
  \end{equation*}
\item Output of $x_T$ for some $T$ large enough.
\end{enumerate}

Here are the correction guarantees.
\begin{thm}[{\citep*{sun_qu_wright}}]\label{thm:sun_qu_wright}
  Let the solution $x^s$ of the phase retrieval problem be arbitrary. Let the measurement vectors $v_1,\dots,v_m$ be generated according to independent normal distributions.

There exists a constant $C>0$ such that, if
\begin{equation*}
Cn\log^3(n)\leq m,
\end{equation*}
then, with probability $1-O\left(\frac{1}{n}\right)$, the sequence of iterates $(x_t)_{t\in\N}$ obtained at the refinement step of the algorithm satisfies
\begin{equation*}
\mathrm{dist}(x^s,x_t)\overset{t\to+\infty}{\to} 0,
\end{equation*}
provided that the gradient descent step is small enough.
\end{thm}

\subsubsection{What are the properties of critical points?}

In order to analyze the algorithm, we first need to understand the properties of critical points. Then we can show that there exists no point with these properties other than the solution $x^s$; in particular, there is no bad critical point.

Let us recall that we call \textit{critical points} the points at which the refinement heuristic (here, gradient descent over $f$) can stagnate. From the definition of gradient descent, a critical point $x_*$ necessarily satisfies
\begin{equation*}
\nabla f(x_*) = 0.
\end{equation*}
This is called a \textit{first-order optimality condition}.

In addition, if we assume that $x_0$ does not belong to some set of ``problematic" initial points which has zero Lebesgue measure, one can show (but it is much more difficult) that critical points must also verify a \textit{second-order optimality condition}:
\begin{equation*}
\nabla^2f(x_*)\succeq 0.
\end{equation*}

These properties, which are not specific to our objective function $f$, are rigorously stated in the following theorem.
\begin{thm}\label{thm:convergence_to_second_order}
Let $\mathcal{L}:\R^n \mbox{ (or $\C^n$) } \to \R$ be analytic. We assume that
\begin{equation*}
\mathcal{L}(x)\to +\infty\quad\mbox{when }||x||\to+\infty.
\end{equation*}
We run gradient descent over $\mathcal{L}$ with constant stepsize $\mu>0$. This yields a sequence of iterates $(x_t)_{t\in\N}$. If $\mu$ is small enough, then, for all $x_0\in B(0,1)$, the sequence $(x_t)_{t\in\N}$ is convergent. Moreover,
\begin{itemize}
\item for all $x_0$, $\nabla \mathcal{L}(\lim_{t\to+\infty} x_t)=0$;
\item for almost all $x_0$, $\nabla^2 \mathcal{L}(\lim_{t\to+\infty} x_t)\succeq 0$.
\end{itemize}
\end{thm}

The first part of this theorem can be deduced from \citep*[Thm 4.1]{absil_mahony_andrews} and the second one from \citep*[Thm 3]{panageas_piliouras}, which is a generalization of \citep*[Cor 9]{lee}. We have stated this theorem for gradient descent with constant stepsize, but it holds for various other local optimization algorithms.

For the algorithm of \citep*{sun_qu_wright}, it implies that:
\begin{itemize}
\item the sequence of iterates $(x_t)_{t\in \N}$ converges to some limite $x_*\in\C^n$ ;
\item with probability 1, $x_*$ satisfies
  \begin{gather}
    \nabla f(x_*) = 0 ; \tag{First-order-optimality} \label{eq:crit_point_first}\\
    \nabla^2 f(x_*) \succeq 0. \tag{Second-order-optimality}.\label{eq:crit_point_second}
  \end{gather}
\end{itemize}

To prove that the algorithm succeeds (Theorem \ref{thm:sun_qu_wright}), it thus suffices to show that, except for the solution $x^s$, no point satisfies Equations \eqref{eq:crit_point_first} and \eqref{eq:crit_point_second}. In the next paragraph, we explain how this can be shown.

\subsubsection{Idea of proof}

To gain some intuition of the proof, let us first determine which points satisfy Equations \eqref{eq:crit_point_first} and \eqref{eq:crit_point_second} when $f$ is replaced with its expectation (which simplifies the computation a log). For any fixed $x\in\C^n$, we can check that the expectation of $f(x)$ over $v_1,\dots,v_m$ is
\begin{equation*}
\E(f(x)) = ||x||^4 - ||x||^2||x^s||^2 - |\scal{x}{x^s}|^2 + ||x^s||^4.
\end{equation*}
Consequently,
\begin{equation*}
\nabla (\E f)(x) = 2((2||x||^2-||x^s||^2)x-\scal{x^s}{x}x^s).
\end{equation*}
With this formula, we can explicitely compute which $x$ satisfy Equation \eqref{eq:crit_point_first} (that is, $\nabla(\E f)(x)=0$). They are the elements of the following three sets:
\begin{itemize}
\item $E_1=\{e^{i\theta} x^s,\theta\in\R\}$, which is the set of solutions of the phase retrieval problem;
\item $E_2=\{0\}$;
\item $E_3=\left\{x\in\C^n,\scal{x^s}{x}=0,||x||=\frac{||x^s||}{\sqrt{2}}\right\}$.
\end{itemize}
Among these points, which ones satisfy Equation \eqref{eq:crit_point_second}? Elements of $E_1$ do, since they are global minimizers of $\E(f)$. For elements of $E_2,E_3$, we need the explicit expression of $\nabla^2(\E f)$: for any $x,h\in\C^n$,
\begin{equation*}
  \nabla^2(\E f)(x)\cdot (h,h) = 2\left((2||x||^2-||x^s||^2)||h||^2+4\mathrm{Re}^2(\scal{x}{h})
  -|\scal{x^s}{h}|^2\right).
\end{equation*}
From this expression, one easily checks that $x=0$ does not satisfy the second-order optimality condition (actually, $\nabla^2(\E f)(0)\prec 0$). And for any $x\in E_3$, we observe that
\begin{equation}
  \nabla^2(\E f)(x)\cdot (x^s,x^s)
  =-2||x^s||^4 < 0.
  \label{eq:hessian_E3_xs}
\end{equation}
As a consequence, Equation \eqref{eq:crit_point_second} cannot hold.

We have thus shown that, replacing $f$ with $\E f$, the only points which satisfy Equations \eqref{eq:crit_point_first} and \eqref{eq:crit_point_second} are the solutions of the phase retrieval problem. In order to extend this result from $\E f$ to $f$ itself, a first idea would be to show that, with high probability, for all $x\in\C^n$,
\begin{subequations}
  \begin{gather}
    ||\nabla f(x)-\nabla\E f(x)||\mbox{ is very small;}\label{eq:sqw_diff_grad} \\
    |||\nabla^2 f(x)-\nabla^2\E f(x)|||\mbox{ is very small,}\label{eq:sqw_diff_hess}
  \end{gather}
\end{subequations}
for some notion of ``smallness'' which would require a proper definition. We could then try to deduce the properties for $f$ from the computations we have done for $\E f$.

This precise approach does not work: whatever way we define ``smallness'', properties \eqref{eq:sqw_diff_grad} and \eqref{eq:sqw_diff_hess} are not true for all $x\in\C^n$ at the same time, with high probability. However, one can still show that $\nabla f$ and $\nabla^2f$ share some properties with their expectations. For instance, with arguments similar to Subsection \ref{ss:no_crit_close}, one can establish the following proposition.
\begin{prop}\label{prop:Z1}
We define, for some explicit constant $\alpha>0$ which is not explicitely given here,
\begin{equation*}
Z_1 = \{x\in\C^n,|\scal{x}{x^s}|\leq\alpha||x^s||^2,||x||\leq (1-\alpha)||x^s||\}.
\end{equation*}

If $m\geq Cn\log^3(n)$ for some constant $C>0$ large enough, it holds with probability $1-O\left(\frac{1}{m}\right)$ that
\begin{equation*}
\nabla^2f(x)\cdot(x^s,x^s)\leq -\alpha ||x^s||^4
\end{equation*}
for all $x\in Z_1$.
\end{prop}
This property, analogous to Equation \eqref{eq:hessian_E3_xs}, for $f$ instead of $\E f$, implies that $f$ has no point in $Z_1$ satisfies Equation \eqref{eq:crit_point_second}. We can define other sets $Z_2,Z_3,\dots$ whose union is equal to $\C^n-E_1$ (recall that $E_1$ is the set of solutions) and establish over $Z_2,Z_3,\dots$ similar properties as the one stated for $Z_1$ in Proposition \ref{prop:Z1}. This proves that $f$ has no point outside $E_1$ satisfies Equations \eqref{eq:crit_point_first} and \eqref{eq:crit_point_second}, and concludes the proof.

\subsubsection{Related work}

Let us emphasize that this proof technique cannot be applied to all algorithms: the non-existence of critical points is a strong property, which many refinement heuristics do not satisfy. In phase retrieval, as far as I know, the algorithm of \citep*{sun_qu_wright} which we have studied is the only one for which the property is known to hold. For alternating projections, for instance, numerical experiments suggest that bad critical points almost always exist\footnote{except if $m\geq O(n^2)$ but this regime is not very interesting}. This does not prevent alternating projections from working well with high probability, for normally distributed measurement vectors, but it makes our proof technique useless.

Outside phase retrieval, this technique has notably been used in \citep*{ge_lee_ma} for matrix completion, in \citep*{bhojanapalli} for RIP matrix sensing, in \citep*{sun_qu_wright_dictionary} for dictionary learning and in \citep*{kawaguchi} for linear neural networks. Many other examples can be found in the review article \citep*{zhang_qu_wright}. Although they have a common proof structure, these results are quite specific to the considered problem and algorithm. Some generalization attempts have been done, notably in \citep*{li_tang} and \citep*{ge_jin_zheng}, but they are still rudimentary.

\section{When there are critical points (leave-one-out)\label{s:leave_one_out}}

In the previous section, we have seen that the success of some non-convex algorithms could be explained by the non-existence of critical points, which removes the stagnation risk for the refinement. We have seen that this phenomenon allows to prove correctness guarantees for various low-rank matrix recovery methods. However, we have also underlined that this proof strategy does not apply to all algorithms, because many succeed \textit{despite the presence of critical points}.

\begin{figure}
  \centering
  \includegraphics[width=0.6\textwidth]{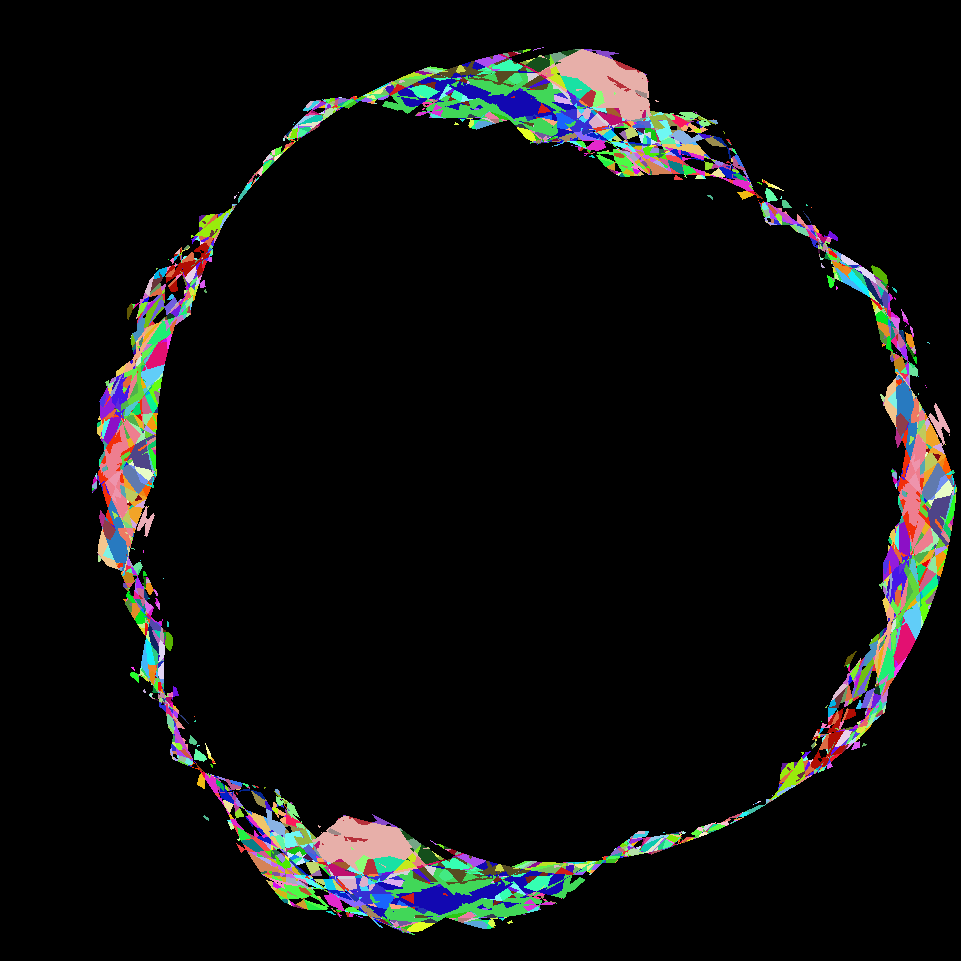}
  \caption{Attraction basins of alternating projections for a phase retrieval problem with $n=20,m=400$
    \label{fig:attraction_basins}}
\end{figure}

An illustration is the alternating projections method for phase retrieval. Figure \ref{fig:attraction_basins} depicts the attraction basins of the various critical points of this algorithm, in a setting where $n=20,m=400$ and the measurement vectors are realizations of independent normal distributions with real coordinates. More precisely, the figure represents a two-dimensional submanifold of $\R^n$, projected onto a square. To each critical point is associated a color; all points in the figure which belong to the attraction basin of this critical point are colored with the corresponding color. The solution of the phase retrieval problem is attributed the color black. Therefore, the black set in Figure \ref{fig:attraction_basins} contains all points starting from which the alternating projections algorithm converges towards the correct solution.

From this figure, we can see that many bad critical points exist. However, the total volume of their attraction basins is tiny: when randomly initialized, alternating projections succeed in finding the correct solution with high probability.

Unfortunately, estimating the size of attraction basins for a given algorithm is in general difficult. If we denote $T$ the operator applied by the refinement heuristic at each iteration, estimating the size of the basins requires to understand the behavior of $(T^n(x_0))_{n\in\N}$ as a function of $x_0$. Exploiting the randomness of $T$ and its independence with $x_0$, it is often possible to understand the properties of $T(x_0)$. But understanding $T^2(x_0)=T(T(x_0))$ (and a fortiori $T^k(x_0)$ for general $k$) is much more difficult: $T$ and $T(x_0)$ are both random variables, but they are correlated and the relation between them is complex.

A proof technique called \textit{leave-one-out} has recently been proposed to overcome this issue. In this section, we present this proof technique through the example of a phase synchronization algorithm, following the article \citep*{zhong_boumal}. At the end, we give a preliminary overview of its possible extensions and limitations, although the leave-one-out technique is too recent so that these are well-understood.

\subsection{Prologue: leave-one-out for machine learning}

This subsection is relatively independent from the rest of the section and can harmlessly be skipped.

\citep*{zhong_boumal} is, as far as I know, the first article in which leave-one-out was used to study a non-convex optimization algorithm. It had however been used in other fields before. In particular, in statistical problems where parameters of a random law must be estimated from samples, it allows to understand some subtle properties of possible estimators \citep*{elkaroui_bean_bickel_lim_yu}. It relies on a relatively simple idea, used for a long time in machine learning, which we briefly describe.

In a classical machine learning problem, the goal is to predict, when given a data point $x\in\R^n$, some property of $x$, modeled by a real number $f_*(x)$. Data are generated according to an unknown probability law $\mathbb{P}$. To learn how to perform the prediction, one is given $m$ ``training'' data points $x_1,\dots,x_m$ (independently generated according to $\mathbb{P}$) and their associated properties $f_*(x_1),\dots,f_*(x_m)$.

Assuming we have designed an algorithm $\mathrm{Alg}$, which takes as input $x_1,\dots$, $x_m$, $f_*(x_1),\dots$, $f_*(x_m)$ and outputs a prediction function $f_{\mathrm{Alg}}:\R^n\to\R$, how do we measure its quality? Ideally, we want $f_{\mathrm{Alg}}$ to precisely mimic $f_*$ at all data points $x$ which could be generated by the probability law $\mathbb{P}$. This motivates the definition of the \textit{generalization error} of $f_{\mathrm{Alg}}$:
\begin{equation*}
\mathrm{Err}(f_{\mathrm{Alg}}) = \E_{x\sim\mathbb{P}}(\ell(f_{\mathrm{Alg}}(x),f_*(x))),
\end{equation*}
for some \textit{loss function} $\ell:\R^2\to\R^+$, well-suited to the problem at hand. Unfortunately, exactly computing the generalization error requires the knowledge of $\mathbb{P}$ and $f_*$. How can we compute an approximate value of this error from only the knowledge of $x_1,\dots,x_m,f_*(x_1),\dots,f_*(x_m)$?

One possibility is to compute, for all $i=1,\dots,m$, the fonction $f_{\mathrm{Alg},-i}$ output by $\mathrm{Alg}$ when $\mathrm{Alg}$ is only given as input$x_1,\dots,x_{i-1},x_{i+1},\dots$, $x_m$ and $f_*(x_1),\dots$, $f_*(x_{i-1})$, $f_*(x_{i+1}),\dots$, $f_*(x_m)$. For any $i$, $f_{\mathrm{Alg},-i}$ should not be very different from $f_{\mathrm{Alg}}$, but it is independent from $(x_i,f_*(x_i))$. One can therefore say that
\begin{equation*}
\ell(f_{\mathrm{Alg},-i}(x_i),f_*(x_i))
\end{equation*}
has almost the same probability distribution as $\ell(f_{\mathrm{Alg}}(x),f_*(x))$ for $x\sim\mathbb{P}$. This leads to the following approximation, called \textit{leave-one-out} approximation:
\begin{equation*}
  \mathrm{Err}(f_{\mathrm{Alg}})\approx \frac{1}{m} \sum_{i=1}^m
  \ell(f_{\mathrm{Alg},-i}(x_i),f_*(x_i)).
\end{equation*}
If the algorithm satisfies some ``stability'' assumptions, this approximation can be rigorously justified \citep*{elisseeff}.

\subsection{Definition of the generalized power method}

In this section, we describe and motivate the phase synchronization algorithm which we will study in the rest of the section.

Let us recall what we have said in Paragraph \ref{sss:phase_sync}: a phase synchronization problem consists in identifying (up to a global phase) $n$ unitary complex numbers $z_1^s,\dots,z_n^s$ from
\begin{equation*}
  C_{k,l}=z_k^s\overline{z_l^s} + w_{k,l},\quad
  \forall k,l\in\{1,\dots,n\},
\end{equation*}
where, for any $k,l$, $w_{k,l}$ is an unknown noise.

We assume to simplify that the $w_{k,l}$ are independent realizations of a complex normal law with variance $\sigma^2$. More precisely, we assume
\begin{align*}
  w_{k,l}&\sim \mathcal{N}_{\C}(0,\sigma^2)
           \mbox{ independently, for all }k<l; \\
  w_{k,l}&=\overline{w_{l,k}}\mbox{ for all }k>l; \\
  w_{k,k}&=0\mbox{ for all }k.
\end{align*}

From the noisy observations $(C_{k,l})_{1\leq k,l\leq n}$ only, it is impossible to exactly recover $(z^s_k)_{1\leq k\leq n}$. Instead, we redefine our goal to finding unitary complex numbers $(z_k^{obj})_{1\leq k\leq n}$ which are minimizers of the following problem:
\begin{equation}\label{eq:def_z_obj}
\min_{|z_1|=\dots=|z_n|=1}\sum_{k,l}|C_{k,l}-z_k\overline{z_l}|^2
\end{equation}
(which corresponds to defining $(z^{obj}_k)_{1\leq k\leq n}$ as the so-called \textit{maximal likelihood} estimator of $(z^s_k)_{1\leq k\leq n}$).

Let us note that Definition \eqref{eq:def_z_obj} is equivalent to $z^{obj}$ being a solution of the following maximization problem:
\begin{equation}\label{eq:def_z_obj_equiv}
z^{obj*}Cz^{obj} = \max_{|z_1|=\dots=|z_n|=1}z^*Cz,
\end{equation}
which is similar to the definition of the main eigenvector of $C$:
\begin{equation}\label{eq:def_z_ppal}
z^{ppal*}Cz^{ppal} = \max_{||z||=1}z^*Cz.
\end{equation}
The main eigenvector can be computed with the \textit{power method}: starting from a random $z^{(0)}$, one iteratively defines $z^{(t+1)}$ as the projection of $Cz^{(t)}$ onto the unit sphere, that is
\begin{equation*}
z^{(t+1)} = \frac{Cz^{(t)}}{||Cz^{(t)}||},\quad\forall t\in\N.
\end{equation*}
The analogy between Problems \eqref{eq:def_z_obj_equiv} and \eqref{eq:def_z_ppal} suggests the following phase synchronization algorithm, similar to the power method.
\begin{enumerate}
\item We set $z^{(0)}=z^{ppal}$.
\item For all $t\geq 0$, we set
  \begin{equation*}
    z^{(t+1)} = \mathcal{P}(Cz^{(t)}),
  \end{equation*}
  where $\mathcal{P}$ is the projection onto $\{z\in\C^n,|z_1|=\dots=|z_n|=1\}$, that is to say, for all $z$,
  \begin{equation*}
    \mathcal{P}(z)_k = \frac{z_k}{|z_k|}, \quad\forall k=1,\dots,n.
  \end{equation*}
  (We use the convention $\frac{0}{0}=1$.)
\end{enumerate}
This algorithm has been introduced in \citep*{boumal} and is named \textit{generalized power method}.

The best known correctness guarantees for this algorithm are stated in the following theorem, which comes from \citep*{zhong_boumal}.
\begin{thm}\label{thm:zhong_boumal}
  There exist constants $\delta>0,\rho\in]0;1[$ such that, if $\sigma \leq \delta \sqrt{\frac{n}{\log(n)}}$, then, with probability $1-O\left(\frac{1}{n^2}\right)$, the sequence of iterates output by the generalized power method satisfies
  \begin{equation*}
    \frac{\mathrm{dist}(z^{(t)},z^{obj})}{||z^{obj}||}\leq \rho^t,\quad
    \forall t\geq 0,
  \end{equation*}
  where $\mathrm{dist}$ is defined by $\mathrm{dist}(u,v)=\inf_{\alpha\in\R}||u-e^{i\alpha}v||$.
\end{thm}

\citep*{gao_zhang} shows that, when $\sigma\geq \sqrt{n}$, $z^{obj}$ is not a significantly better approximation of $z^s$ than a random guess: for some constant $c>0$,
\begin{equation*}
\E\left(\mathrm{dist}(z^{obj},z^s)^2\right)\geq cn,
\end{equation*}
while a random vector $\hat z$ chosen with uniform probability in $\{z,|z_1|=\dots=|z_n|=1\}$ satisfies $\E\left(\mathrm{dist}(\hat z,z^s)^2\right)\leq 2n$. This means that, when $\sigma\geq \sqrt{n}$, our phase synchronization problem is rather uninteresting. Consequently, the condition $\sigma\leq \delta \sqrt{\frac{n}{\log(n)}}$ in Theorem \ref{thm:zhong_boumal} is not very restrictive.

\subsection{Proof of Theorem \ref{thm:zhong_boumal}}

\subsubsection{Difficulties\label{sss:difficulties}}
  
The goal of this paragraph is to explain why proving Theorem \ref{thm:zhong_boumal} is difficult and requires a different strategy from Section \ref{s:no_bad_critical}. Pointing out the main difficulties will also help us understand the intuition which led Zhong and Boumal to a correct proof, presented in Paragraph \ref{sss:proof}. Towards this goal, let us try to apply the strategy of Section \ref{s:no_bad_critical} (that is to show that the generalized power method has no critical point) and see why it fails.

We start with a basic but necessary property (whose proof is in Appendix \ref{s:fixed_point}).
\begin{prop}\label{prop:fixed_point}
The vector $z^{obj}$ is a fixed point of the operator $(z\to \mathcal{P}(Cz))$.
\end{prop}

To follow the strategy of Section \ref{s:no_bad_critical}, let us try to prove that $(z\to\mathcal{P}(Cz))$ has no bad critical point in some neighborhood of $z^s$. As neighborhood, let us take the simplest choice $B(z^s,\epsilon||z^s||)$, for some $\epsilon>0$ which can be much smaller than $1$ but must not go to $0$ when $n$ goes to infinity. We could try to establish the following properties:
\begin{enumerate}
\item The operator $(z\to \mathcal{P}(Cz))$ is $\rho$-Lipschitz on $B(z^s,\epsilon||z^s||)$ for some $\rho<1$: for all $y,y'\in B(z^s,\epsilon||z^s||)$,
  \begin{equation}\label{eq:PC_lipschitz}
    \mathrm{dist}(\mathcal{P}(Cy),\mathcal{P}(Cy')) \leq \rho \mathrm{dist}(y,y').
  \end{equation}
\item Vectors $z^{obj},z^0$ belong to $B\left(z^s,\frac{\epsilon}{3}||z^s||\right)$.
\end{enumerate}
Combined with Proposition \ref{prop:fixed_point}, these two properties would imply that, for all $t$,
\begin{equation*}
\mathrm{dist}(z^{(t)},z^{obj})\leq \rho^t \mathrm{dist}(z^{(0)},z^{obj}),
\end{equation*}
and thus prove the result.

The second property turns out to be true. Unfortunately, the first one is false with high probability. Indeed, $\mathcal{P}$ is discontinuous at any point $z\in\C^n$ with at least one zero coordinate. If we want Equation \eqref{eq:PC_lipschitz} to be true for all $y,y'\in B(z^s,\epsilon||z^s||)$, it is therefore necessary that, for any $y\in B(z^s,\epsilon ||z^s||)$, no coordinate of $Cy$ is equal to zero. Can this happen?

Let us first consider a vector $y=z^s+u$, for some random vector $u$ chosen in $B(0,\epsilon ||z^s||)$ with uniform probability, independently from the noise $(w_{k,l})_{1\leq k,l\leq n}$. We denote $W$ the matrix whose $(k,l)$-th coefficient is $w_{k,l}$. Let us recall that it is a Hermitian matrix and that its off-diagonal coefficients are realizations of centered Gaussian variables with variance $\sigma^2$. The matrix $C$ is
\begin{equation*}
C = z^sz^{s*}+W.
\end{equation*}
For all $k$, we thus have
\begin{align}
  (C(z^s+u))_k
  & = z_k^s\left(||z^s||^2+\scal{z^s}{u}\right) + \scal{W_{:,k}}{z^s+u}\nonumber \\
  & = z_k^s\left(||z^s||^2 + O(\epsilon||z^s||^2)\right) + O(\sigma||z^s+u||)\nonumber \\
  & = nz_k^s + O(\epsilon n) + O(\sigma \sqrt{n}).\label{eq:C_z_plus_u}
\end{align}
In these equations, $W_{:,k}$ is the $k$-th column of $W$. The second equality is due to the fact that $||u||\leq\epsilon||z^s||$ on the one hand and, on the other hand, to the fact that, since $u$ and $W$ are independent, $\scal{W_{:,k}}{z^s+u}$ is, conditioning on $u$, a centered Gaussian variable with variance $||z^s+u||^2\sigma^2$; its values are therefore of order $\sigma ||z^s+u||$. The last equality is true because $||z^s||=\sqrt{n}$.

As $|z_k^s|=1$ and $\sigma\ll \sqrt{n}$, Equation \eqref{eq:C_z_plus_u} yields
\begin{equation*}
|C(z^s+u)|_k = n(1+O(\epsilon)),
\end{equation*}
which is far away from $0$. This informal reasoning suggests that, if we pick a vector $y\in B(z^s,\epsilon||z^s||)$ at random with uniform probability, independently from $W$, then $Cy$ has no coordinate equal or close to zero, with high probability.

Nevertheless, if we consider $y=z^s-\eta \frac{W{:,k}}{||W_{:,k}||}$ (for some complex number $\eta$ such that $|\eta|<\epsilon||z^s||$ whose value will be given later),
\begin{align*}
  (Cy)_k
  & = z^s_k\left(||z^s||^2 - \eta \scal{z^s}{\frac{W_{:,k}}{||W_{:,k}||}}\right)
    + \scal{W_{:,k}}{z^s} - \eta ||W_{:,k}|| \\
  & = z^s_k n (1+O(\epsilon)) - \eta ||W_{:,k}|| \\
  & = z^s_k n(1+O(\epsilon)) - \eta \sigma \sqrt{n} (1+o(1)).
\end{align*}
We choose $\eta = \frac{\sqrt{n}(1+O(\epsilon))}{\sigma(1+o(1))}z_k^s$. This leads to
\begin{equation*}
  (Cy)_k=0.
\end{equation*}
As $||z^s||=\sqrt{n}$, this definition of $\eta$ satisfies the constraint $|\eta|<\epsilon||z^s||$ as soon as $\sigma > \frac{2}{\epsilon}$, so $y$ does belong to $B(z^s,\epsilon||z^s||)$ and one of the coordinates of $Cy$ is zero. And since the theorem must be proved for $\sigma$ of order up to $\sqrt{\frac{n}{\log(n)}}$, we cannot assume $\sigma \leq \frac{2}{\epsilon}$.

To summarize, Property \eqref{eq:PC_lipschitz} is plausible for vectors $y,y'$ chosen at random in $B(z^s,\epsilon||z^s||)$ independently from $W$, but false when $y$ or $y'$ is unnaturally correlated with a column of $W$.

\subsubsection{Proof\label{sss:proof}}

Hiding some technical aspects under the carpet, we can say that the proof of Theorem \ref{thm:zhong_boumal} decomposes in the following four steps:
\begin{enumerate}
\item For some $\rho\in]0;1[$, we show that the map $(z\to\mathcal{P}(Cz))$ is $\rho$-Lipschitz over
  \begin{equation}\label{eq:def_N}
    \mathcal{N}\overset{def}{=}\{z\in B(z^s,\epsilon||z^s||)\mbox{ such that }\forall k,|\scal{W_{:,k}}{z}|<\kappa n\}.
  \end{equation}
  (In this definition, $\epsilon,\kappa>0$ are well-chosen constants.)
\item We show that, with high probability, $z^{(0)}$ belongs to
  \begin{equation*}
    \mathcal{N}'\overset{def}{=}\left\{z\in B\left(z^s,\frac{\epsilon}{3}||z^s||\right)\mbox{ such that }\forall k,|\scal{W_{:,k}}{z}|<\frac{\kappa n}{2}\right\}
    \subset\mathcal{N}.
  \end{equation*}
\item We show that, with high probability, $z^{(t)}$ belongs to $\mathcal{N}'$ for all $t\in\N^*$.
\item We conclude: from the previous three steps, $(z^{(t)})_{t\in\N}$ is a Cauchy sequence. It is therefore convergent, and its limit belongs to $\mathcal{N}$. To establish the theorem, we only have to prove that the limit is $z^{obj}$.
\end{enumerate}

The first of these steps is the least difficult one: $\mathcal{N}$ is precisely defined as (more or less) the largest set over which the classical probabilistic arguments informally used in Paragraph \ref{sss:difficulties} allow to show that $(z\to \mathcal{P}(Cz))$ is $\rho$-Lipschitz for some $\rho\in]0;1[$. The first step therefore only requires to make rigorous the arguments of Paragraph \ref{sss:difficulties}. The fourth step relies on a relation between the definition of $z^{obj}$ and some semidefinite optimization problem. We will not describe it in more detail. The second and third steps are the ones for which the leave-one-out technique is necessary. As the main objective of this section is to present the leave-one-out principle and not to provide a complete proof of the theorem, we assume that the property of the second step is true and focus on proving the property of the third step.

We proceed iteratively. Assuming that $z^{(0)},\dots,z^{(t-1)}$ belong to $\mathcal{N}'$, we must show that $z^{(t)}$ belongs to $\mathcal{N}'$. We can relatively easily show that
\begin{equation*}
z^{(t)}\in B\left(z^s,\frac{\epsilon}{3}||z^s||\right).
\end{equation*}
The key difficulty is to show that, for all $k\leq n$,
\begin{equation}\label{eq:scal_Wk_z}
  |\scal{W_{:,k}}{z^{(t)}}|<\frac{\kappa n}{2}.
\end{equation}
Following the reasoning of Paragraph \ref{sss:difficulties}, it would be easy if, for any $k$, $z^{(t)}$ and $W_{:,k}$ were independent random variables. However, since the construction of $z^{(t)}$ involves the matrix $C$, $z^{(t)}$ depends on $C$, and therefore on $W$.

To overcome this issue, we introduce, for each $k\leq n$, an auxiliary sequence $(z^{(k,t)})_{t\in\N}$, whose definition is exactly the same as $(z^{(t)})_{t\in\N}$, except that we replace the matrix
\begin{equation*}
C = z^sz^{s*}+W
\end{equation*}
with the matrix
\begin{equation*}
C^{(k)}=z^sz^{s*}+W^{(k)},
\end{equation*}
where $W^{(k)}$ is the same matrix as $W$, with the coefficients in the $k$-th row and column replaced by zeroes. Let us underline that these auxiliary sequences are theoretical tools only: they cannot be computed in practice. Indeed, in a real phase synchronization instance, computing an auxiliary sequence $(z^{(k,t)})_{t\in\N}$ would require computing $W^{(k)}$ and thus exactly knowing $W$ (which is never the case, otherwise the problem would be trivial).

For all $k$, the sequence $(z^{(k,t)})_{k\in\N}$ is independent from $W_{:,k}$, since the coefficients in $W_{:,k}$ appear nowhere in the definition of the sequence. As a consequence, it holds for all $t$, with high probability, that
\begin{equation}\label{eq:scal_Wk_zk}
|\scal{W_{:,k}}{z^{(k,t)}}| = \tilde{O}(\sigma\sqrt{n})<\frac{\kappa n}{4}.
\end{equation}
(Compared to $O$, the notation $\tilde O$ hides an additional $\sqrt{\log(n)}$ factor.)

To deduce Property \eqref{eq:scal_Wk_z} from Property \eqref{eq:scal_Wk_zk}, it suffices to show that, for any $k$, sequences $(z^{(t)})_{t\in\N}$ and $(z^{(k,t)})_{t\in\N}$ are close with high probability. More precisely, we prove by iteration over $t$ that the following properties are true with high probability\footnote{For technical reasons, Zhong and Boumal proceed by iteration only up to $t\approx 3n^2$. For larger values of $t$, they use a different, more elementary, argument, which we do not describe here.}:
\begin{enumerate}
\item $\mathrm{dist}(z^{(t)},z^{(k,t)})\leq \frac{1}{60}$ for all $k$;\label{it:loo_small_dist}
\item $|\scal{W_{:,k}}{z^{(t)}}|<\frac{\kappa n}{2}$ for all $k$.\label{it:loo_scal_prod}
\end{enumerate}

Let us give an idea of the proof of these properties. For all $k$,
\begin{align}
  \mathrm{dist}(z^{(t)},z^{(k,t)})
  & = \mathrm{dist}(\mathcal{P}(C z^{(t-1)}),\mathcal{P}(C^{(k)}z^{(k,t-1)})) \nonumber \\
  & \leq \mathrm{dist}(\mathcal{P}(C z^{(t-1)}),\mathcal{P}(C z^{(k,t-1)}))
    + \mathrm{dist}(\mathcal{P}(C z^{(k,t-1)}),\mathcal{P}(C^{(k)} z^{(k,t-1)})).
    \label{eq:triangular_inequality}
\end{align}
Using our induction hypotheses,
\begin{equation*}
  z^{(t-1)}\in\mathcal{N}'
  \quad\mbox{and}\quad
  \mathrm{dist}(z^{(k,t-1)},z^{(t-1)})\leq \frac{1}{60},
\end{equation*}
from which we deduce that $z^{(t-1)},z^{(k,t-1)}$ belong to $\mathcal{N}$ and therefore that
\begin{equation}
  \mathrm{dist}(\mathcal{P}(Cz^{(t-1)}),\mathcal{P}(Cz^{(k,t-1)}))
  \leq \rho\mathrm{dist}(z^{(t-1)},z^{(k,t-1)})\leq \frac{\rho}{60}.
  \label{eq:first_term}
\end{equation}
On the other hand, with high probability
\begin{align}
  \mathrm{dist}(\mathcal{P}(C z^{(k,t-1)}),\mathcal{P}(C^{(k)} z^{(k,t-1)}))
  & \overset{(*)}{\leq} \frac{2}{n} \mathrm{dist}(C z^{(k,t-1)},C^{(k)} z^{(k,t-1)}) \nonumber\\
  & \leq \frac{2}{n} ||(C-C^{(k)})z^{(k,t-1)} || \nonumber\\
  & \overset{(\square)}{\leq} \frac{2}{n} \left( |\scal{W_{:,k}}{z^{(k,t-1)}}| + ||W_{:,k}||\,|z^{(k,t-1)}_k|
    \right) \nonumber\\
  & \overset{(\circ)}{=} \tilde{O} \left(\frac{\sigma}{\sqrt{n}}\right).
    \label{eq:second_term}
\end{align}
Inequality $(*)$ is true because the coordinates of $Cz^{(k,t-1)}$ and $C^{(k)}z^{(k,t-1)}$ are all larger than $\frac{n}{2}$ in modulus (which can be derived from the fact that $z^{(t-1)}$ and $z^{(k,t-1)}$ belong to $\mathcal{N}$). Inequality $(\square)$ is true because all rows and columns of $C-C^{(k)}$ are zero, except the $k$-th column (which is $W_{:,k}$) and the $k$-th row (which is $W^*_{:,k}$). For inequality $(\circ)$, we have used the independence between $W_{:,k}$ and $z^{(k,t-1)}$ to upper bound $|\scal{W_{:,k}}{z^{(k,t-1)}}|$.

Combining inequalities \eqref{eq:triangular_inequality}, \eqref{eq:first_term} and \eqref{eq:second_term} yields
\begin{equation*}
  \mathrm{dist}(z^{(t)},z^{(k,t)})
  \leq \frac{\rho}{60} + \tilde{O}\left(\frac{\sigma}{\sqrt{n}}\right)
  \leq \frac{1}{60}.
\end{equation*}

This proves the Property \ref{it:loo_small_dist}. Property \ref{it:loo_scal_prod} is easier: for all $k$,
\begin{align*}
  |\scal{W_{:,k}}{z^{(t)}}|
  & \leq |\scal{W_{:,k}}{z^{(k,t)}}| + ||W_{:,k}|| \mathrm{dist}(z^{(t)},z^{(k,t)}) \\
  & = \tilde{O} (\sigma \sqrt{n}) + \tilde{O}(\sigma\sqrt{n}) \times \frac{1}{60} \\
  & < \frac{\kappa n}{2}.
\end{align*}

\subsection{Extensions and limits}

This subsection briefly describes how leave-one-out can be applied to other non-convex algorithms than the generalized power method and when its application fails. The reader must keep in mind that, as said before, the introduction of leave-one-out in the field of non-convex optimization is recent. Therefore, the applicability and limits of this method have not been well explored yet, and the content of this subsection must be considered as preliminary ideas.

\subsubsection{General principle}

To understand how leave-one-out can be applied to other algorithms, we briefly summarize its principle in a general context. We consider a non-convex algorithm, whose sequence of iterates $(z^{(t)})_{t\in\N}$ is defined by
\begin{equation*}
z^{(t+1)}=T(z^{(t)})
\end{equation*}
for some map $T:\R^n\to\R^n$. We want to show that it converges to the ``correct solution" $z^{obj}$.

The principle of leave-one-out is as follows.
\begin{enumerate}
\item We define an appropriate set, of the form
  \begin{equation*}
    \mathcal{N}=\{z\in\Omega\mbox{ such that }F_1(z)<1,\dots,F_K(z)<1\},
  \end{equation*}
  for some functions $F_1,\dots,F_K:\R^n\to\R$ and $\Omega$ a (deterministic) open set. In the case of the generalized power method, this corresponds to Definition \eqref{eq:def_N}, with $F_k:z\to\frac{1}{\kappa n}|\scal{W_{:,k}}{z}|$ for $k=1,\dots,n$ and $\Omega=B(z^s,\epsilon||z^s||)$.
\\
We show that
\begin{itemize}
\item $T$ is contractive over $\mathcal{N}$ (or some other property which implies that $z^{(t)}\overset{t\to+\infty}{\longrightarrow} z^{obj}$ if $z^{(t)}$ belongs to $\mathcal{N}$ for all $t$);
\item $F_1,\dots,F_K$ are Lipschitz in a neighborhood of any point of $\mathcal{N}$.
\end{itemize}
\item We show that $z^{(0)}\in\mathcal{N}$.
\item We introduce auxiliary iterates $(z^{(k,t)})_{1\leq k\leq K,t\in\N}$, obeying a similar definition as $(z^{(t)})_{t\in\N}$ except that $T$ is replaced by variants $T_k$: for any $k,t$,
\begin{equation*}
z^{(k,t+1)}=T_k(z^{(k,t)}).
\end{equation*}
The definition of these auxiliary sequences must be chosen in such a way that, for any $k,t$, it is not difficult to upper bound
\begin{equation*}
F_k(z^{(k,t)}).
\end{equation*}
\item We prove by iteration over $t\in\N$ that, for well-chosen real numbers $(\epsilon_{k,t})_{1\leq k\leq K}$,
  \begin{gather*}
    \forall k,\quad ||z^{(t)}-z^{(k,t)}|| \leq \epsilon_{k,t}; \\
    z^{(t)}\in\mathcal{N}.
  \end{gather*}
  For the first of these properties, the principle of the proof is broadly the following sequence of inequalities:
  \begin{align*}
    ||z^{(t)}-z^{(k,t)}||
    & = ||T(z^{(t-1)})-T_k(z^{(k,t-1)})|| \\
    & \leq ||T(z^{(t-1)})-T(z^{(k,t-1)})|| + ||T(z^{(k,t-1)})-T_k(z^{(k,t-1)})|| \\
    & \leq \rho_T \epsilon_{k,t-1} + ||(T-T_k)(z^{(k,t-1)})||,
  \end{align*}
  where $\rho_T$ is the Lipschitz constant of $T$ over $\mathcal{N}$. As $T_k$ is a variant of $T$, one can hope to be able to upper bound $||(T-T_k)(z^{(k,t-1)})||$ by something ``small".
  \\
  For the second property, the difficulty is to show that $F_k(z^{(t)})<1$ for all $k$. We have
  \begin{align*}
    F_k(z^{(t)})
    & \leq F_k(z^{(k,t)}) + \left|F_k(z^{(t)})-F_k(z^{(k,t)})\right| \\
    & \leq F_k(z^{(k,t)}) + \rho_{F_k} \epsilon_{k,t},
  \end{align*}
  where $\rho_{F_k}$ is the Lipschitz constant of $F_k$ in the neighborhood of $z^{(t)}$. From the well-chosen definition of $(z^{(k,t)})_{t\in\N}$, we have a good upper bound for $F_k(z^{(k,t)})$ at our disposal, which allows to conclude.
\end{enumerate}

This is the general principle applied, for instance, in \citep*{ma_wang_chi_chen,ding_chen,chen_chi_fan_ma,chen_chi_fan_ma_yan}. It can be further refined. For instance, in \citep*{chen_chi_fan_ma}, the map $T$ is not exactly contractive, which makes it difficult to upper bound $||T(z^{(t-1)})-T(z^{(k,t-1)})||$; the authors achieve this through a second leave-one-out argument, nested in the first one.

\citep*{zhang_alt_min} can also be seen as a (quite extreme) instance of this principle, applied to alternating projections for phase retrieval in a setting where the number of measurement vectors is much larger than the dimension of the unknown signal. In this article, $T_k$ is not a variant of $T$: it is a random map, with the same distribution as $T$ but independent from it.

\subsubsection{Limits}

A first drawback of this method is that it requires a precise understanding of the local Lipschitz properties of $T$, possibly also of $dT$. This makes the proof rather technical, and very specific to one problem and one algorithm.

In my opinion, a second drawback is that it strongly relies on the fact that $T$ is continuous, and even contractive, on a relatively large set $\mathcal{N}$. This is necessary if want to guarantee that sequences $(z^{(t)})_{t\in\N}$ and $(z^{(k,t)})_{t\in\N}$ stay close to each other. But some algorithms do not satisfy this contraction property.

\begin{figure}
  \centering
    \begin{tikzpicture}
    \begin{axis}[
      xmin = 0,
      ymin = 0, ymax = 0.15,
      xticklabels = {?,0,0.02,0.04,0.06,0.08,0.1},
      yticklabels = {?,0,0.05,0.1,0.15},
      xlabel = {$||z-z'||$},
      ylabel = {$||T(z)-T(z')||$},
      legend pos = outer north east]

      \addplot+ coordinates {
        (0.0025, 0.01629128576296882)
        (0.01, 0.03715796368147183)
        (0.025, 0.06147254838861682)
        (0.05, 0.09006631239191809)
        (0.075, 0.1133533293884338)
        (0.1, 0.13431554428213968)
      };

      \addplot+ coordinates {
        (0.0025, 0.0017039957974030098)
        (0.01, 0.006807011776764489)
        (0.025, 0.017044994488547686)
        (0.05, 0.03406297986696872)
        (0.075, 0.05108538780929064)
        (0.1, 0.06818228612076621)
      };
      
      \legend{Alternating projections, Wirtinger Flow};
    \end{axis}
  \end{tikzpicture}
\caption{$\E ||T(z)-T(z')||$ as a function of $||z-z'||$ for two phase retrieval algorithms, \textit{Wirtinger Flow} and alternating projections, for the reconstruction of signals with dimension $n=400$ from $m=4000$ phaseless measurements. For any distance $||z-z'||$, $\E ||T(z)-T(z')||$ is computed by averaging over $1000$ random pairs $(z,z')$. \label{fig:limits}}
\end{figure}
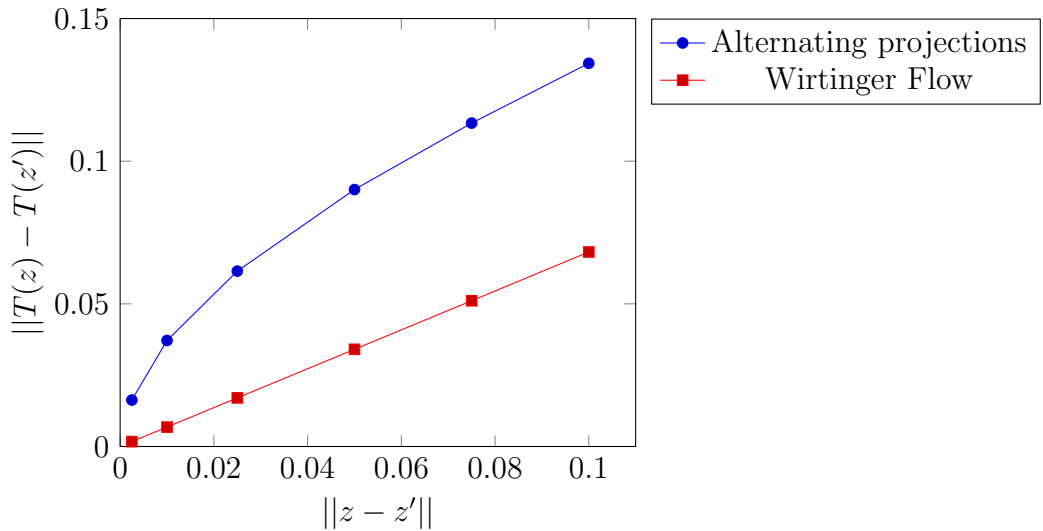

An example is (again) the alternating projections method for phase retrieval, applied to signals and measurement vectors with real coordinates\footnote{We emphasize the realness. It seems that the behaviour is quite different when coordinates are complex.}, which follow independent normal laws. Indeed, a rough computation suggests that, for this algorithm, for any points $z,z'$ on the unit sphere,
\begin{equation*}
\mathbb{E}(||T(z)-T(z')||^2)
\mbox{ is of order }\frac{n}{m}||z-z'||,
\end{equation*}
which suggests that $T$ is far from being locally Lipschitz on a relatively large set (otherwise the expectation should be of order $||z-z'||^2$). Figure \ref{fig:limits} supports this assertion: it shows that, for the \textit{Wirtinger Flow} algorithm described in Subsection \ref{ss:no_crit_close} (to which leave-one-out can be applied \citep*{chen_chi_fan_ma}), $||T(z)-T(z')||$ is essentially proportional to $||z-z'||$ and the proportionality constant is smaller than $1$. For alternating projections, $||T(z)-T(z')||$ grows much faster at small values of $||z-z'||$.

\section{Burer-Monteiro methods\label{s:burer_monteiro}}

In this final section, we study a family of non-convex algorithms, the so-called \textit{Burer-Monteiro methods}, introduced in \citep*{burer_monteiro03}. They are applicable to all low-rank matrix recovery problems satisfying reasonably general assumptions. Compared to the previous two sections, where we described techniques to study algorithms applied to specific low-rank problems satisfying restrictive statistical assumptions, the correction guarantees we present in this section have a much larger application field.

\subsection{Definition of Burer-Monteiro methods}

In this subsection, we first describe the class of low-rank matrix recovery problems to which Burer-Monteiro methods are applicable (Paragraph \ref{sss:problems}). We then define the Burer-Monteiro methods themselves (Paragraph \ref{sss:burer_monteiro_def}) and give a first overview of their behavior through a basic numerical experiment (Paragraph \ref{sss:burer_monteiro_num}).

\subsubsection{Problems\label{sss:problems}}

As previously, we consider problems where one must recover a low-rank matrix $X^s$ from simple information, modeled by the fact that $X^s$ belongs to some set $\mathcal{E}$:
\begin{equation*}
  \mbox{minimize }\mathrm{rank}(X)\mbox{ for }X\in\mathcal{E}.
  \tag{\ref{eq:min_rank}}
\end{equation*}

Burer-Monteiro methods can be seen as a ``deconvexification" of the convexified techniques presented in Section \ref{s:convex}. Therefore, they only apply to problems which admit a \textit{convex relaxation} satisfying precise properties. The first of these properties is that the convex relaxation must be \textit{exact}: Problem \eqref{eq:min_rank} and its convex approximation must have the same solution $X^s$, so that solving the convex version suffices to solve Problem \eqref{eq:min_rank}.

\begin{hyp}\label{hyp:exact}
Problem \eqref{eq:min_rank} has an exact convex relaxation.
\end{hyp}

The second property is as follows.

\begin{hyp}\label{hyp:SDP}
The convex relaxation can be written under the form
\begin{equation}\label{eq:relaxation}
\mbox{minimize }\Tr(CX)\mbox{ for }X\in\mathcal{E}_{SDP},
\end{equation}
for some $C\in\mathcal{S}_n(\R)$ (called the \textrm{cost matrix}) and $\mathcal{E}_{SDP}$ a subset of $\mathcal{S}_n(\R)$ which is compact and is defined as the intersection between $\mathcal{S}_n^+(\R)$ and an affine space of dimension $m\geq 1$.
\end{hyp}

When this assumption is satisfied, Problem \eqref{eq:relaxation} is a semidefinite program:
  \begin{align}
    \mbox{minimize }&\Tr(CX), \nonumber\\
    \mbox{with }&\mathcal{A}(X)=b, \tag{SDP}\label{eq:SDP}\\
    &X\succeq 0,\nonumber
  \end{align}
  for some linear map $\mathcal{A}:\mathcal{S}_n(\R)\to\R^m$ and some vector $b\in\R^m$ such that
\begin{equation*}
  \mathcal{E}_{SDP}=\{X\in\mathcal{S}_n(\R),\mathcal{A}(X)=b,X\succeq 0\}.
\end{equation*}

Let us note that, in this section, we constrain all matrices which come into play to have real coefficients. This restriction simply aims at simplifying the proofs: up to minor modifications, our main results also hold true for complex coefficients.

All low-rank recovery problems described in the introduction satisfy Assumptions \ref{hyp:exact} and \ref{hyp:SDP}.
\begin{enumerate}
\item We have seen (Equation \eqref{eq:convex_completion}) that matrix completion problems could be approximated by the following convex problem:
  \begin{align*}
  \mbox{minimize }&||X||_*\\
  \mbox{with }&X_{i,j} = X^s_{ij},\quad \forall (i,j)\in\Omega.
  \end{align*}
  This problem can be reformulated under the form \eqref{eq:SDP}\footnote{This guarantees that Assumption \ref{hyp:SDP} holds, except possibly for the compactness requirement (but it turns out to be unnecessary in this case).} because of the following equality, which holds for all $X\in\R^{n_1\times n_2}$:
  \begin{equation*}
    ||X||_* = \min\left\{
      \frac{\Tr(Y)+\Tr(Z)}{2},Y\in\mathcal{S}_{n_1}(\R),Z\in\mathcal{S}_{n_2}(\R),
      \left(\begin{smallmatrix}Y&X\\X^T&Z\end{smallmatrix}\right)\succeq 0
    \right\}.
  \end{equation*}
  In addition, we have seen in Subsection \ref{ss:convex_completion} that, under appropriate assumptions, the convex relaxation is exact with high probability.
\item We have already seen that phase retrieval problems admit several convex relaxations of the form \eqref{eq:SDP}:
  \begin{center}
    \vskip -0.5cm
    \begin{minipage}{0.45\textwidth}
      \begin{align*}
        & \mbox{\eqref{eq:PhaseLift}} \\
      \mbox{minimize }&\Tr(X), \\
        \mbox{with }& \scal{X}{v_kv_k^*} = |\scal{x^s}{v_k}|^2, \forall k\leq m, \\
                    & X\succeq 0.
      \end{align*}
    \end{minipage}\quad
    \begin{minipage}{0.45\textwidth}
      \begin{align*}
        & \mbox{\eqref{eq:PhaseCut}} \\
        \mbox{minimize }&\Tr(MU), \\
        \mbox{with }& U_{k,k}=1, \forall k\leq m,\\
        &U\succeq 0.
      \end{align*}
    \end{minipage}
  \end{center}
  We have seen that, at least in a specific statistical setting, \eqref{eq:PhaseLift} is exact with high probability. This is also true for \eqref{eq:PhaseCut}.
\item Phase synchronization problems admit a convex relaxation of the form
  \begin{align*}
    \mbox{minimize }&-\Tr(CU), \\
    \mbox{with }&U_{k,k}=1,\quad\forall k\leq m, \\
    &U\succeq 0.
  \end{align*}
  In the case where phase measures are contaminated with a Gaussian additive noise (which is the setting we have studied in Section \ref{s:leave_one_out}), the relaxation is exact with high probability, provided that the noise level satisfies\footnote{We observe that this assumption is the same as in Theorem \ref{thm:zhong_boumal}. This is not a coincidence: the fact that the non-convex generalized power method succeeds in recovering $z^{obj}$ is a crucial tool for establishing the exactness of the convex relaxation.} $\sigma\leq c\sqrt{\frac{n}{\log(n)}}$ for some constant $c>0$ \citep*{zhong_boumal}.
\end{enumerate}

\subsubsection{Definition of Burer-Monteiro methods\label{sss:burer_monteiro_def}}

Let us assume that Assumptions \ref{hyp:exact} and \ref{hyp:SDP} hold and discuss (again) the question of how to solve Problem \eqref{eq:SDP}:
\begin{gather*}
  \mbox{minimize }\Tr(CX)\mbox{ for }X\in\mathcal{E}_{SDP},
  \tag{\ref{eq:SDP}} \\
  \mbox{where }\mathcal{E}_{SDP}=\{X\in\mathcal{S}_n(\R),\mathcal{A}(X)=b,X\succeq 0\}.\nonumber
\end{gather*}
As mentioned in Section \ref{s:convex}, many general solvers with rigorous convergence guarantees exist for problems of this form, but they tend to be too slow for medium to high-dimensional applications.

Burer-Monteiro methods, on the other hand, are of heuristic nature. They sometimes fail at correctly identifying the desired minimizer but, when they succeed, they can provide significant speed-ups.
Their principle is to ``deconvexify" Problem \eqref{eq:SDP}. At first sight, this can appear counter-intuitive: since Problem \eqref{eq:SDP} has been constructed by ``convexifying" a low-rank matrix recovery problem, isn't there a risk that, when deconvexifying it, we fall back to a non-convex problem essentially equivalent to the initial one? Actually, no. After the deconvexification, we arrive at a family of non-convex problems, parametrized by an integer $p$, and some members of this family can be easier to solve than the initial formulation.

The ``deconvexification" relies on the observation that, since the convex relaxation is exact (Assumption \ref{hyp:exact}), the solution $X^s$ of Problem \eqref{eq:SDP} has low rank. For all $p\in\N^*$, we define
\begin{equation*}
\mathcal{E}_{SDP,p}=\mathcal{E}_{SDP}\cap \{X\in\mathcal{S}_n(\R),\mathrm{rank}(X)\leq p\}.
\end{equation*}
For any $p\geq \mathrm{rank}(X^s)$, Problem \eqref{eq:SDP} is equivalent to
\begin{equation*}
\mbox{minimize }\Tr(CX)\mbox{ for }X\in\mathcal{E}_{SDP,p}.
\end{equation*}
Informally, if $p\ll n$, the set $\mathcal{E}_{SDP,p}$ has a much smaller dimension than $\mathcal{E}_{SDP}$. We can therefore hope that minimizing our cost function $X\to \Tr(CX)$ over $\mathcal{E}_{SDP,p}$ requires significantly less computational effort than minimizing it over $\mathcal{E}_{SDP}$.

To take advantage of this dimensionality reduction, the simplest idea is to parametrize $\mathcal{E}_{SDP,p}$ by a low-dimensional manifold $\mathcal{M}_p$. Specifically, we consider the following mapping:
\begin{equation*}
V\in\mathcal{M}_p \quad\to \quad V V^T \in\mathcal{E}_{SDP,p},
\end{equation*}
with
\begin{equation*}
\mathcal{M}_p=\{V\in\R^{n\times p},\mathcal{A}(VV^T)=b\}.
\end{equation*}
This mapping is onto: any element $X\in\mathcal{E}_{SDP,p}$ can be factorized as $X=VV^T$ for some $V\in\R^{n\times p}$, since it is a positive semidefinite matrix with rank at most $p$. And the equality $\mathcal{A}(X)=b$ implies that $\mathcal{A}(VV^T)=b$ and thus that $V$ belongs to $\mathcal{M}_p$.

Using this parametrization, we can rewrite Problem \eqref{eq:SDP} as
\begin{equation}\label{eq:SDP_fact}
  \mbox{minimize }f_C(V)
  \overset{def}{=} \Tr(CVV^T)\mbox{ for }V\in\mathcal{M}_p.
  \tag{Factorized SDP}
\end{equation}
In this version of the problem, the unknown $V$ has $np$ coefficients, which is much less than the $n^2$ coefficients of the original unknown $X$ if $p\ll n$. Manipulating $V$ is therefore less costly than manipulating $X$. However, compared to \eqref{eq:SDP}, Problem \eqref{eq:SDP_fact} has the major drawback of not being convex anymore. As a consequence, bad critical points may exist and simple algorithms applicable to Problem \eqref{eq:SDP_fact} are not guaranteed to succeed, although it is possible that they work beautifully.

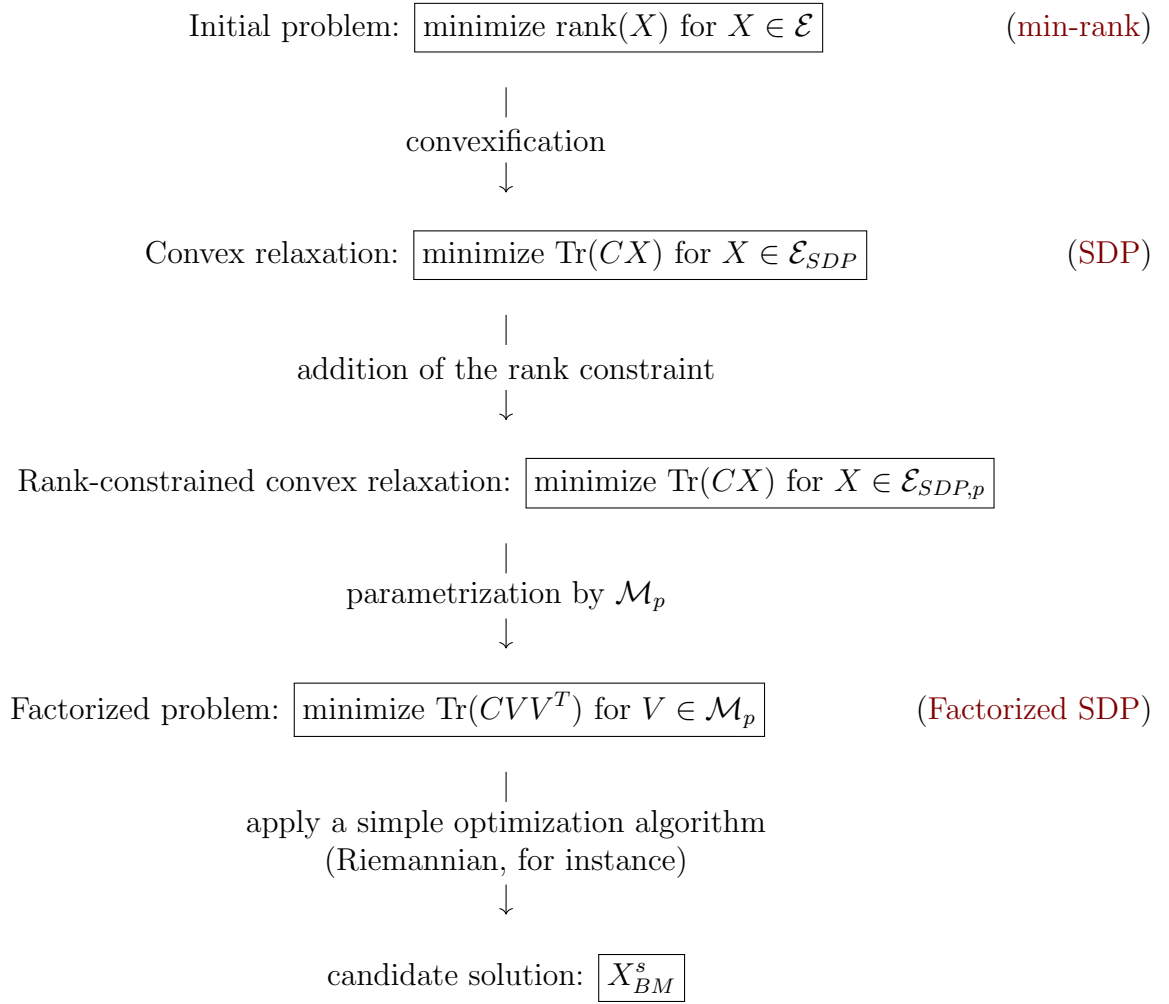
\begin{figure}
  \centering
  \begin{equation*}
    \mbox{Initial problem: }
    \boxed{\mbox{minimize }\mathrm{rank}(X)\mbox{ for }X\in\mathcal{E}}
    \tag{\ref{eq:min_rank}}
  \end{equation*}
  \begin{tikzpicture}
    \draw[->] (0,-0.3) -- (0,-0.7) ;
    \draw (0,0.3) -- (0,0.7) ;
    \node at (0,0) {convexification};
  \end{tikzpicture}
  \begin{equation*}
    \mbox{Convex relaxation: }
    \boxed{\mbox{minimize }\Tr(CX)\mbox{ for }X\in\mathcal{E}_{SDP}}
    \tag{\ref{eq:SDP}}
  \end{equation*}
  \begin{tikzpicture}
    \draw[->] (0,-0.3) -- (0,-0.7) ;
    \draw (0,0.3) -- (0,0.7) ;
    \node at (0,0) {addition of the rank constraint};
  \end{tikzpicture}
  \begin{equation*}
    \mbox{Rank-constrained convex relaxation: }
    \boxed{\mbox{minimize }\Tr(CX)\mbox{ for }X\in\mathcal{E}_{SDP,p}}
  \end{equation*}
  \begin{tikzpicture}
    \draw[->] (0,-0.3) -- (0,-0.7) ;
    \draw (0,0.3) -- (0,0.7) ;
    \node at (0,0) {parametrization by $\mathcal{M}_p$};
  \end{tikzpicture}
  \begin{equation*}
    \mbox{Factorized problem: }
    \boxed{\mbox{minimize }\Tr(CVV^T)\mbox{ for }V\in\mathcal{M}_p}
    \tag{\ref{eq:SDP_fact}}
  \end{equation*}
  \begin{tikzpicture}
    \draw[->] (0,-0.8) -- (0,-1.2) ;
    \draw (0,0.3) -- (0,0.7) ;
    \node at (0,0) {apply a simple optimization algorithm};
    \node at (0,-0.5) {(Riemannian, for instance)};
  \end{tikzpicture}
  \begin{equation*}
    \mbox{candidate solution: }\boxed{X^s_{BM}}
  \end{equation*}
\caption{Schematic view of a Burer-Monteiro method\label{fig:burer_monteiro}}
\end{figure}

We call \textit{Burer-Monteiro method} any solver which attempts to solve Problem \eqref{eq:min_rank} by applying any reasonable algorithm to the factorized problem \eqref{eq:SDP_fact}. This algorithmic scheme is summarized in Figure \ref{fig:burer_monteiro}. In these notes, we limit ourselves to Burer-Monteiro methods where Problem \eqref{eq:SDP_fact} is solved by \textit{Riemannian optimization}. This requires the set $\mathcal{M}_p$ to be a Riemannian submanifold\footnote{A reader unfamiliar with Riemannian geometry can imagine a Riemannian submanifold as a regular curve or surface inside $\R^{n\times p}$.} of $\R^{n\times p}$. We actually need a slightly stronger assumption.

\begin{hyp}\label{hyp:submersion}
  For all $V_0\in\mathcal{M}_p$, the differential at $V_0$ of the mapping
  \begin{equation*}
    \tilde{\mathcal{A}}:V \in\R^{n\times p} \quad \to\quad
    \mathcal{A}(VV^T)\in\R^m
  \end{equation*}
  is onto.
\end{hyp}

\begin{prop}
If Assumption \ref{hyp:submersion} is satisfied, $\mathcal{M}_p$ is a Riemannian submanifold of $\R^{n\times p}$, with dimension $np-m$.
\end{prop}

Going back to the examples mentioned at the end of Paragraph \ref{sss:problems}, the last two among them (that is phase retrieval, for the \eqref{eq:PhaseCut} formulation, and phase synchronization) satisfy this assumption. For matrix completion, it is less clear; it probably depends on the ground truth $X^s$.

Riemannian optimization algorithms are in general local optimization methods, which start at some point of the considered manifold, and progressively ``move'', using the gradient (and possibly the Hessian) of the cost function to decide the movement direction. They oftentimes derive from a classical optimization algorithm over $\R^d$. For instance, two of the most prominent Riemannian algorithms are \textit{Riemannian gradient descent} (which derives from gradient descent over $\R^d$) and \textit{Riemannian Trust-Region method} (which derives from the Trust-Region method over $\R^d$). Many others exist \citep*{absil}. Each one of them can be more or less adapted to a given application. Here, we will try to keep our discussion general, and will not make any particular assumption on the Riemannian optimization algorithm.

\subsubsection{Numerical experiment\label{sss:burer_monteiro_num}}

Our goal, in the rest of this section, is to provide a partial answer to the following question:
\begin{center}
  When can we guarantee that a Burer-Monteiro method succeeds?
\end{center}
We focus on guarantees which require essentially no hypothesis besides Assumptions \ref{hyp:exact}, \ref{hyp:SDP} and \ref{hyp:submersion}, so that our results apply to as many \eqref{eq:min_rank} problems as possible. In particular, $C,\mathcal{A},b$ can be arbitrary.

Nevertheless, we cannot avoid making an hypothesis on $p$. We must choose this hypothesis with great care, since it determines the practical relevance of the correctness guarantees. Indeed, in Problem \eqref{eq:SDP_fact}, the unknown $V$ has dimension proportional to $p$ so Burer-Monteiro methods run much faster when $p$ is small. As far as I know, large values of $p$ are rarely used in practice, so it is important for correctness guarantees to hold for as small values of $p$ as possible. On the other hand, whether Problem \eqref{eq:SDP_fact} has bad critical points or not may strongly depend on $p$. Consequently, it may be that establishing guarantees for small values of $p$ is very difficult or even impossible.

\begin{figure}
  \begin{subfigure}{0.48\textwidth}
    \centering
    \resizebox{\textwidth}{!}{
    \begin{tikzpicture}
      \begin{axis}[
        xmin = 0,
        ymin = -0.01, ymax = 1.01,
        xlabel = {$m/n$},
        ylabel = {Success rate},
        legend pos = north west]

        \addplot+ coordinates { (0.0, 0.0) (1.0, 0.0) (2.0, 0.0) (3.0,
          0.0) (4.0, 0.4) (5.0, 1.0) (6.0, 1.0) (7.0, 1.0) (8.0,
          1.0) };
          
        \addplot+ coordinates { (0.0, 0.0) (1.0, 0.0) (2.0, 0.0) (3.0,
          0.4) (4.0, 0.44) (5.0, 0.8) (6.0, 0.95) (7.0, 1.0) (8.0,
          1.0) };

        \addplot+ coordinates { (0.0, 0.0) (1.0, 0.0) (2.0, 0.0) (3.0,
          0.0) (4.0, 1.0) (5.0, 1.0) (6.0, 1.0) (7.0, 1.0) (8.0, 1.0)
        };

        \legend{\eqref{eq:PhaseCut},$p=1$,$p=2$};
      \end{axis}
    \end{tikzpicture}}
  \end{subfigure}
  \begin{subfigure}{0.48\textwidth}
    \centering
    \resizebox{\textwidth}{!}{
    \begin{tikzpicture}
      \begin{axis}[
        xmin = 0,
        ymin = -0.01, ymax = 1.01,
        xlabel = {$m/n$},
        ylabel = {Success rate},
        legend pos = north west]
          
        \addplot+ coordinates { (0.0, 0.0) (1.0, 0.0) (2.0, 0.0) (3.0,
          0) (4.0, 0.25) (5.0, 1) (6.0, 1) (7.0, 1) (8.0, 1) };

        \addplot+ coordinates { (0.0, 0.0) (1.0, 0.0) (2.0, 0.0) (3.0,
          0.0) (4.0, 0.05) (5.0, 0.) (6.0, 0.05) (7.0, 0.05) (8.0,
          0.25) };
      
        \addplot+ coordinates { (0.0, 0.0) (1.0, 0.) (2.0, 0.0) (3.0,
          0.15) (4.0, 0.95) (5.0, 1) (6.0, 1.0) (7.0, 1) (8.0,
          1.0) };
        
        \legend{\eqref{eq:PhaseCut},$p=1$,$p=2$};
      \end{axis}
    \end{tikzpicture}}
  \end{subfigure}
  \caption{Success rate, as a function of $m/n$, for three algorithms: \eqref{eq:PhaseCut}, a Burer-Monteiro method with $p=1$, a Burer-Monteiro method with $p=2$. The signal dimension is $n=32$. The figure on the left corresponds to measurement vectors independently chosen according to normal laws and the figure on the right to measurement vectors representing a ``wavelet transform''. Success rates have been computed by averaging the results of $20$ reconstruction attempts.
    \label{fig:bm_phasecut}}
\end{figure}
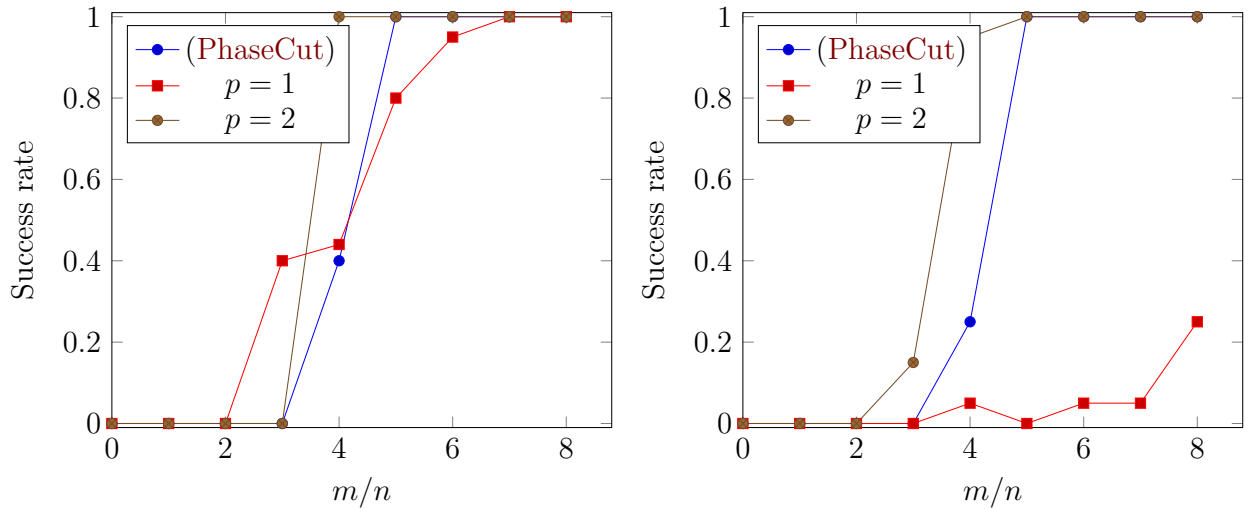

To get a preliminary intuition, let us study the numerical performance of a Burer-Monteiro method on toy phase retrieval problems. The specific method we have chosen for this experiment uses the \eqref{eq:PhaseCut} relaxation of the phase retrieval problem, and solves the corresponding \eqref{eq:SDP_fact} by Riemannian gradient descent. In this case, when the relaxation is exact, its solution has rank $1$, so $p$ can a priori take any positive integer value. In our tests, we use either $p=1$ or $p=2$. Signals have dimension $n=32$ and are chosen according to a complex normal distribution. We vary the number of measurements between $m=0$ and $m=8n$. Results are displayed on Figure \ref{fig:bm_phasecut}. As a reference point, the figure also displays the success rate of the convex \eqref{eq:PhaseCut} method (which consists in directly solving Problem \eqref{eq:SDP} with an interior-point method, without considering Problem \eqref{eq:SDP_fact}).

The experiment covers two types of measurement vectors: realizations of independent normal vectors on the one hand and vectors describing a ``wavelet transform'' on the other hand. Vectors of the second type are more ``structured'' (in particular, some of them are very close to being orthogonal to each other while others are strongly correlated). It is known that phase retrieval problems with one or the other measurement type have significantly different intrinsic properties (regarding stability to noise, for instance); problems of the second type generally tend to be more difficult. Hence, it is possible that Burer-Monteiro methods do not behave the same for the two measurement types; this is the reason for including both in our experiment.

From Figure \ref{fig:bm_phasecut}, we see that, when $p=1$, the behavior of our Burer-Monteiro method indeed depends on the measurement type: it works fine (that is, essentially as well as \eqref{eq:PhaseCut}) for normal vectors, and fails for the wavelet transform. On the contrary, when $p=2$, it works as well as \eqref{eq:PhaseCut} in both cases.

Much more extensive numerical experiments can be found in the literature, involving different Burer-Monteiro methods and other low-rank recovery problems than phase retrieval \citep*{burer_monteiro03,journee,boumal_riemannian}. Overall, they lead to a conclusion compatible with our toy experiment: some failure cases can be observed if $p$ is equal or almost equal to the rank of the solution of Problem \eqref{eq:SDP}. However, it suffices for $p$ to be ``slightly larger'' than the rank for these failure cases to disappear.

\subsection{Success guarantees when $\frac{p(p+1)}{2}>m$}

The observations of Paragraph \ref{sss:burer_monteiro_num} suggest that Burer-Monteiro have excellent empirical success rates, even for quite small values of $p$. Can we prove it?

In this subsection, we describe the correctness guarantees established in \citep*{boumal_voroninski_bandeira}, which can informally be stated as follows: under Assumptions \ref{hyp:exact}, \ref{hyp:SDP}, \ref{hyp:submersion}, provided that
\begin{equation*}
  \frac{p(p+1)}{2}>m,
\end{equation*}
Riemannian algorithms succeed at solving almost any problem of the form \eqref{eq:SDP_fact} (that is, Burer-Monteiro methods succeed). We recall that $m$ is the number of affine constraints in Problem \eqref{eq:SDP_fact}, that is, the dimension of the range of $\mathcal{A}$.

In the following paragraph, we precisely state these guarantees and in the next ones, we give an idea of how to prove them.

\subsubsection{Precise statement}

We first recall a result of Section \ref{s:no_bad_critical}, Theorem \ref{thm:convergence_to_second_order}: if we run gradient descent over an analytic function $f:\R^d\to\R$, we obtain (under relatively weak assumptions) a converging sequence of iterates, whose limit $x_*$ is a second-order critical point of $f$:
\begin{equation*}
  \nabla f(x_*) = 0
  \quad\mbox{and}\quad
  \nabla^2f(x_*)\succeq 0.
\end{equation*}
It turns out that this property is also true for Riemannian algorithms, at least part of them \citep*{boumal_absil_cartis,criscitiello_boumal}: when applied to Problem \eqref{eq:SDP_fact}, these algorithms are guaranteed to find a matrix $V_*$ which is a second-order critical point of $f_C$:
\begin{equation*}
  \nabla_{\mathcal{M}_p}f_C(V_*)=0
  \quad\mbox{and}\quad
  \nabla^2_{\mathcal{M}_p}f_C(V_*)\succeq 0.
\end{equation*}
Here, $\nabla_{\mathcal{M}_p}$ and $\nabla^2_{\mathcal{M}_p}$ are respectively the gradient and Hessian of $f_C$, restricted to the manifold $\mathcal{M}_p$.

Therefore, we can establish correctness guarantees for Burer-Monteiro methods with the same argument as in Subsection \ref{ss:no_crit_at_all}: the minimizers\footnote{I use plural since Problem \eqref{eq:SDP_fact} never has a single minimizer: if $V_*$ is a minimizer, so is $V_* G$ for any orthogonal $p\times p$ matrix $G$.} of Problem \eqref{eq:SDP_fact} are second-order critical points of $f_C$ (it is a general property of minimizers). If they are the \textit{only} second-order critical points, Burer-Monteiro methods are guaranteed to correctly solve Problem \eqref{eq:SDP_fact} and thus the initial problem \eqref{eq:min_rank}.

\begin{thm}[{\citep*{boumal_voroninski_bandeira}}]\label{thm:bvb}
  Let $\mathcal{A}:\mathcal{S}_n(\R)\to\R^m,b\in\R^m$ be fixed. We assume that Assumptions \ref{hyp:exact}, \ref{hyp:SDP} and \ref{hyp:submersion} are verified and that
  \begin{equation*}
    \frac{p(p+1)}{2}>m.
  \end{equation*}
  Then, for all cost matrices $C\in\mathcal{S}_n(\R)$ outside some zero Lebesgue measure set, Problem \eqref{eq:SDP_fact} has no second-order critical point except its minimizers.
\end{thm}

The following paragraphs contain an overview of the proof of this theorem. Parts of it are relatively technical and possibly difficult to understand, but I hope that they allow, even without being read in full detail, to get an idea of which tools are necessary to prove this result, and why the condition $\frac{p(p+1)}{2}>m$ appears.

\subsubsection{Proof principle for Theorem \ref{thm:bvb}}

Let us assume that $\mathcal{A},b$ are fixed and Assumptions \ref{hyp:exact}, \ref{hyp:SDP} and \ref{hyp:submersion} are verified. Theorem \ref{thm:bvb} follows from the following two lemmas.
\begin{lem}\label{lem:rank_deficient_crit_point}
  For all $V_*\in\mathcal{M}_p$, whatever the cost matrix $C$, if
  \begin{equation*}
    \nabla_{\mathcal{M}_p}f_C(V_*)=0,\quad
    \nabla^2_{\mathcal{M}_p}f_C(V_*)\succeq 0\quad
    \mbox{and}\quad
    \mathrm{rank}(V_*)<p,
  \end{equation*}
  then $V_*$ is a minimizer of Problem \eqref{eq:SDP_fact}.
\end{lem}

\begin{lem}\label{lem:no_crit_point_rank_p}
  We assume that $\frac{p(p+1)}{2}>m$.

  For any matrix $C\in\mathcal{S}_n(\R)$ outside some zero Lebesgue measure set, there does not exist $V_*\in\mathcal{M}_p$ such that
  \begin{equation*}
    \nabla_{\mathcal{M}_p}f_C(V_*)=0\quad
    \mbox{and}\quad
    \mathrm{rank}(V_*)=p.
  \end{equation*}
\end{lem}

Indeed, if these two lemmas are true, it means that, when $\frac{p(p+1)}{2}>m$, for all cost matrices outside some zero Lebesgue measure set:
\begin{itemize}
\item Problem \eqref{eq:SDP_fact} has no second-order critical point (nor even a first-order one) with rank $p$ (Lemma \ref{lem:no_crit_point_rank_p});
\item Problem \eqref{eq:SDP_fact} has no second-order critical point with rank strictly smaller than $p$, except its solutions (Lemma \ref{lem:rank_deficient_crit_point}).
\end{itemize}
Consequently, Problem \eqref{eq:SDP_fact} has no second-order critical point except its solutions (let us note that elements of $\mathcal{M}_p$ have $p$ columns and are therefore of rank at most $p$).

We will not explain the proof of Lemma \ref{lem:rank_deficient_crit_point}. We simply present a rough geometrical interpretation for it. Let us imagine that the second-order critical points of $f_C$ are exactly its local minimizers (it is not true in full generality, but this hypothesis helps developing an intuition). If $V_*$ is a local minimizer of $f_C$ over $\mathcal{M}_p$, then $V_*V_*^T$ is a local minimizer of $X\to \scal{C}{X}$ over $\mathcal{E}_{SDP,p}$ (since $V \in\mathcal{M}_p \to VV^T\in\mathcal{E}_{SDP,p}$ is a parametrization of $\mathcal{E}_{SDP,p}$).

Now, Lemma \ref{lem:rank_deficient_crit_point} is true because of the following (non obvious) property: if $\mathrm{rank}(V_*V_*^T)<p$, the set of possible ``displacement directions'' around $V_*V_*^T$ in $\mathcal{E}_{SDP,p}$ forms a cone of elements of $\mathcal{S}_n(\R)$, whose convex envelope contains $\mathcal{E}_{SDP}$. When $V_*V_*^T$ is a local minimizer of $X\to\scal{C}{X}$ over $\mathcal{E}_{SDP,p}$, we have $\scal{C}{X-V_*V_*^T}\geq 0$ for all $X$ in the ``displacement cone'', and thus also for all $X$ in the convex envelope of the cone. Therefore, if the envelope contains $\mathcal{E}_{SDP}$, we have $\scal{C}{X-V_*V_*^T}\geq 0$ for all $X\in\mathcal{E}_{SDP}$. In particular, $V_*V_*^T$ is a solution of Problem \eqref{eq:SDP}; this implies that $V_*$ is a solution of Problem \eqref{eq:SDP_fact}.

\subsubsection{Idea of proof for Lemma \ref{lem:no_crit_point_rank_p}}

We want to show that, when $\frac{p(p+1)}{2}>m$, there exists no $V_*\in\mathcal{M}_p$ such that
\begin{equation}\label{eq:crit_point_rank_p}
  \nabla_{\mathcal{M}_p}f_C(V_*)=0\quad
  \mbox{and}\quad
  \mathrm{rank}(V_*)=p,
\end{equation}
except possibly if $C$ belongs to some zero Lebesgue measure set.

Let us first give an explicit description of the set of cost matrices $C$ for which, on the contrary, there exists $V_*$ satisfying Equation \eqref{eq:crit_point_rank_p}. Then we will show that this set has Lebesgue measure zero.

For any full rank $V_*\in\mathcal{M}_p$, the map $V\in\mathcal{M}_p\to VV^T\in\mathcal{E}_{SDP,p}$ is essentially\footnote{It is not exactly a diffeomorphism since the map is invariant to right multiplication with any orthogonal matrix. To get an exact diffeomorphism, we must quotient $\mathcal{M}_p$ by the set of orthogonal matrices.} a diffeomorphism in some neighborhood of $V_*$ (because one can check that the rank of its differential is locally constant). The equality $\nabla_{\mathcal{M}_p}f_C(V_*)=0$ is therefore equivalent to the fact that the gradient of $X\in\mathcal{E}_{SDP,p}\to\scal{C}{X}$ is zero at $V_*V_*^T$. In other words, it is equivalent to $C$ being orthogonal to the tangent space to $\mathcal{E}_{SDP,p}$ at $V_*V_*^T$.

Consequently, the set of cost matrices $C$ for which there exists $V_*$ satisfying Equation \eqref{eq:crit_point_rank_p} is included in (actually, equal to)
\begin{equation}\label{eq:union_orthogonals}
\bigcup_{M\in\mathcal{E}_{SDP,p},\mathrm{rank}(M)=p}(T_M\mathcal{E}_{SDP,p})^\perp.
\end{equation}
(For any $M\in\mathcal{E}_{SDP,p}$ with rank $p$, we denote by $T_M\mathcal{E}_{SDP,p}$ the tangent space to $\mathcal{E}_{SDP,p}$ at $M$.)

If $N$ is the dimension of $\mathcal{S}_n(\R)$ (that is, $\frac{n(n+1)}{2}$) and $D$ the dimension of the manifold $\{M\in\mathcal{E}_{SDP,p},\mathrm{rank}(M)=p\}$, the vector space $(T_M\mathcal{E}_{SDP,p})^{\perp}$ has dimension $N-D$ for any $M\in\mathcal{E}_{SDP,p}$ with rank $p$. As a consequence, the set in Equation \eqref{eq:union_orthogonals} is the union, parametrized by a $D$-dimensional manifold, of $(N-D)$-dimensional spaces. Intuitively, as illustrated by Figure \ref{fig:non_degenerated_union}, it is a set with ``dimension''\footnote{I use quotes because the set is not a manifold and therefore has no ``dimension'' in the usual sense of this word.} at most
\begin{equation*}
(N-D)+D = N.
\end{equation*}

We have shown that the set of ``problematic'' cost matrices is an at most $N$-dimensional subset of the $N$-dimensional vector space $\mathcal{S}_n(\R)$. An $N$-dimensional subset of an $N$-dimensional space has no reason to have zero Lebesgue measure: this conclusion does not seem very useful. Fortunately, it is possible to refine the computation of the dimension with the help of the following proposition.
\begin{prop}\label{prop:segment}
Let us assume that $\frac{p(p+1)}{2}>m$. For any $M\in\mathcal{E}_{SDP,p}$ with rank $p$, there exists a segment in $\mathcal{E}_{SDP,p}$ such that $M$ is in the interior of the segment and the tangent space to $\mathcal{E}_{SDP,p}$ is constant on the segment.
\end{prop}

This proposition implies that, in Equation \eqref{eq:union_orthogonals}, the parametrization by rank $p$ elements of $\mathcal{E}_{SDP,p}$ is redundant. The set can actually be written as a union of $(N-D)$-dimensional vector spaces, parametrized by a set with dimension only $D-1$. The ``dimension'' of this set is therefore at most
\begin{equation*}
(N-D) + (D-1) = N-1,
\end{equation*}
which means that it is a zero Lebesgue measure subset of $\mathcal{S}_n(\R)$. Figure \ref{fig:degenerated_union} illustrates this argument.

\begin{figure}[h]
  \begin{subfigure}{0.47\textwidth}
    \centering
    \begin{tikzpicture}
      \draw[domain=0:400,smooth, variable=\th] plot ({3*cos(\th)},{2*sin(\th)}) ;
      \draw[domain=0:180,dashed,smooth, variable=\th] plot ({3*cos(\th)},{0.4*sin(\th)}) ;
      \draw[domain=180:360,smooth, variable=\th] plot ({3*cos(\th)},{0.4*sin(\th)}) ;
      \draw[-{Latex[length=3mm]},red] ({3/sqrt(3)},{2/sqrt(3)}) -- ({13/3/sqrt(3)},{4/sqrt(3)});
      \draw[-{Latex[length=3mm]},red] (0,2) -- (0,4);
      \draw[-{Latex[length=3mm]},red] (1,-1) -- (15/9,-2.5);
      \draw[-{Latex[length=3mm]},red] ({-3/sqrt(2)},{2/sqrt(2)}) -- ({-13/3/sqrt(2)},{4/sqrt(2)});
      \draw[-{Latex[length=3mm]},red] ({-3/sqrt(2)},{-0.4/sqrt(2)}) -- ({-6/sqrt(2)},{-1.3/sqrt(2)});
    \end{tikzpicture}
    \caption{A $2$-dimensional manifold included in $\R^3$. The union of lines which are orthogonal to one of its tangent spaces is the whole ambiant space $\R^3$; it has dimension $3$.
    \label{fig:non_degenerated_union}}
  \end{subfigure}\quad
  \begin{subfigure}{0.47\textwidth}
    \begin{tikzpicture}
      \draw[blue] ({-3*cos(55)+0.6*sin(55)*0.286},{-3*sin(55)*sqrt(1.04)*0.958}) --
      ({5*cos(55)-sin(55)*0.286},{5*sin(55)*sqrt(1.04)*0.958});
      \draw[blue] ({-3*cos(55)+0.6*sin(55)*0.286},{3*sin(55)*sqrt(1.04)*0.958}) --
      ({5*cos(55)-sin(55)*0.286},{-5*sin(55)*sqrt(1.04)*0.958});
      \draw[domain=0:180,blue,smooth,dashed,variable=\th] plot ({2*cos(55)+0.4*sin(55)*sin(\th)},
      {2*sqrt(1.04)*sin(55)*cos(\th)});
      \draw[domain=180:360,blue,smooth,variable=\th] plot ({2*cos(55)+0.4*sin(55)*sin(\th)},
      {2*sqrt(1.04)*sin(55)*cos(\th)});
      \draw[domain=0:180,blue,smooth,dashed,variable=\th] plot ({-2*cos(55)+0.4*sin(55)*sin(\th)},
      {-2*sqrt(1.04)*sin(55)*cos(\th)});
      \draw[domain=180:360,blue,smooth,variable=\th] plot ({-2*cos(55)+0.4*sin(55)*sin(\th)},
      {-2*sqrt(1.04)*sin(55)*cos(\th)});
      \draw (0,0) -- ({6*cos(35)-1.2*sin(35)*0.14},{6*sin(35)*sqrt(1.04)*0.99});
      \draw (0,0) -- ({6*cos(35)-1.2*sin(35)*0.14},{-6*sin(35)*sqrt(1.04)*0.99});
      \draw[domain=0:180,smooth,dashed,variable=\th] plot ({3*cos(35)+0.6*sin(35)*sin(\th)},
      {3*sqrt(1.04)*sin(35)*cos(\th)});
      \draw[domain=180:360,smooth,variable=\th] plot ({3*cos(35)+0.6*sin(35)*sin(\th)},
      {3*sqrt(1.04)*sin(35)*cos(\th)});
      \draw[-{Latex[length=3mm]},red] (
      {5*cos(35)-sin(35)*0.14},{5*sin(35)*sqrt(1.04)*0.99}) --
      ({5*cos(35)-sin(35)*0.14-1.5*sin(35)-0.3*cos(35)*0.14},
      {5*sin(35)*sqrt(1.04)*0.99+1.5*cos(35)*sqrt(1.04)*0.958});
      \draw[-{Latex[length=3mm]},red] (
      {3.5*cos(35)-0.7*sin(35)*0.14},{3.5*sin(35)*sqrt(1.04)*0.99}) --
      ({3.5*cos(35)-0.7*sin(35)*0.14-1.5*sin(35)-0.3*cos(35)*0.14},
      {3.5*sin(35)*sqrt(1.04)*0.99+1.5*cos(35)*sqrt(1.04)*0.958});
      \draw[-{Latex[length=3mm]},red] (
      {2*cos(35)-0.4*sin(35)*0.14},{2*sin(35)*sqrt(1.04)*0.99}) --
      ({2*cos(35)-0.4*sin(35)*0.14-1.5*sin(35)-0.3*cos(35)*0.14},
      {2*sin(35)*sqrt(1.04)*0.99+1.5*cos(35)*sqrt(1.04)*0.958});      
      \draw[-{Latex[length=3mm]},red] (
      {5*cos(35)-sin(35)*0.14},{-5*sin(35)*sqrt(1.04)*0.99}) --
      ({5*cos(35)-sin(35)*0.14-1.5*sin(35)-0.3*cos(35)*0.14},
      {-5*sin(35)*sqrt(1.04)*0.99-1.5*cos(35)*sqrt(1.04)*0.958});
      \draw[-{Latex[length=3mm]},red] (
      {3.5*cos(35)-0.7*sin(35)*0.14},{-3.5*sin(35)*sqrt(1.04)*0.99}) --
      ({3.5*cos(35)-0.7*sin(35)*0.14-1.5*sin(35)-0.3*cos(35)*0.14},
      {-3.5*sin(35)*sqrt(1.04)*0.99-1.5*cos(35)*sqrt(1.04)*0.958});
      \draw[-{Latex[length=3mm]},red] (
      {2*cos(35)-0.4*sin(35)*0.14},{-2*sin(35)*sqrt(1.04)*0.99}) --
      ({2*cos(35)-0.4*sin(35)*0.14-1.5*sin(35)-0.3*cos(35)*0.14},
      {-2*sin(35)*sqrt(1.04)*0.99-1.5*cos(35)*sqrt(1.04)*0.958});      
      \draw[ultra thin] (0,0) -- ({6*cos(35)-1.2*sin(35)*cos(18)},
      {-6*sqrt(1.04)*sin(35)*sin(18)});
      \draw[-{Latex[length=3mm]},red] ({2*cos(35)-0.4*sin(35)*cos(18)},
      {-2*sqrt(1.04)*sin(35)*sin(18)}) -- ({2*cos(35)-0.4*sin(35)*cos(18)-1.5*sin(35)
      -1.5*0.2*cos(35)*cos(18)},{-2*sqrt(1.04)*sin(35)*sin(18)-1.5*sqrt(1.04)*sin(18)*cos(35)});
      \draw[-{Latex[length=3mm]},red] ({3.5*cos(35)-0.7*sin(35)*cos(18)},
      {-3.5*sqrt(1.04)*sin(35)*sin(18)}) -- ({3.5*cos(35)-0.7*sin(35)*cos(18)-1.5*sin(35)
      -1.5*0.2*cos(35)*cos(18)},{-3.5*sqrt(1.04)*sin(35)*sin(18)-1.5*sqrt(1.04)*sin(18)*cos(35)});
      \draw[-{Latex[length=3mm]},red] ({5*cos(35)-sin(35)*cos(18)},
      {-5*sqrt(1.04)*sin(35)*sin(18)}) -- ({5*cos(35)-sin(35)*cos(18)-1.5*sin(35)
      -1.5*0.2*cos(35)*cos(18)},{-5*sqrt(1.04)*sin(35)*sin(18)-1.5*sqrt(1.04)*sin(18)*cos(35)});
    \end{tikzpicture}
    \centering
    \caption{A manifold (here a cone without its vertex, in black) whose points are all inside some segment along which the tangent space is constant. The union of lines which are orthogonal to one of its tangent spaces has dimension $2$ (it is the cone drawn in blue).\label{fig:degenerated_union}}
  \end{subfigure}
\end{figure}

\subsubsection{Idea of proof for Proposition \ref{prop:segment}}

Let $M\in\mathcal{E}_{SDP,p}$ with rank $p$ be fixed.

We want to show that there is a segment in $\mathcal{E}_{SDP,p}$, containing $M$, over which the tangent space to $\mathcal{E}_{SDP,p}$ does not vary. For this, we need an explicit expression for the tangent space.

Let us recall that
\begin{equation*}
  \mathcal{E}_{SDP,p} = \{X\succeq 0\}
  \cap\{X,\mathcal{A}(X)=b\}\cap\{X,\mathrm{rank}(X)\leq p\}.
\end{equation*}
Since $M$ is positive semidefinite and $\mathrm{rank}(M)=p$, all symmetric matrices close enough to $M$ have rank at least $p$ and all rank $p$ symmetric matrices close enough to $M$ are positive semidefinite. Consequently, in some neighborhood of $M$, $\mathcal{E}_{SDP,p}$ coincides with
\begin{equation*}
\{X,\mathcal{A}(X)=b\}\cap\{X,\mathrm{rank}(X)=p\}.
\end{equation*}
From this remark (and Assumption \ref{hyp:submersion}), we can show that
\begin{align*}
  T_M\mathcal{E}_{SDP,p}
  & = \mathrm{Ker}(\mathcal{A})\cap T_M\{X,\mathrm{rank}(X)=p\} \\
  & = \mathrm{Ker}(\mathcal{A})\cap \{X\in\mathcal{S}_n(\R),\scal{v}{Xv}=0\mbox{ for all }v\in\mathrm{Range}(M)^\perp\}.
\end{align*}
We see that the tangent space depends on $M$ only through the range of $M$. Therefore, to show that the tangent space is constant on a segment, it suffices to show that there exists a segment of $\mathcal{M}_p$ (containing $M$) whose elements all have the same range.

The set $\{X\in\mathcal{S}_n(\R),\mathrm{Im}(X)\subset\mathrm{Im}(M)\}$ has dimension
\begin{equation*}
\frac{\dim(M)(\dim(M)+1)}{2} = \frac{p(p+1)}{2}.
\end{equation*}
Since $\frac{p(p+1)}{2}>m=\dim(\mathrm{Range}(\mathcal{A}))$, the set must contain a matrix $H\ne 0$ such that
\begin{equation*}
\mathcal{A}(H)=0.
\end{equation*}
We fix such a matrix $H$.

For any $t\in\R$ close enough to $0$,
\begin{itemize}
\item $\mathrm{rank}(M+tH)\geq \mathrm{rank}(M)=p$ and $\mathrm{Range}(M+tH)\subset\mathrm{Range}(M)$, from which we deduce
  \begin{equation*}
    \mathrm{rank}(M+tH)=p \quad\mbox{and}\quad
    \mathrm{Range}(M+tH)=\mathrm{Range}(M);
  \end{equation*}
\item $M+tH\succeq 0$;
\item $\mathcal{A}(M+tH)=\mathcal{A}(M)=b$.
\end{itemize}
In particular, for $\epsilon>0$ small enough, $[M-\epsilon H;M+\epsilon H]$ is a segment included in $\mathcal{E}_{SDP,p}$ whose elements all have the same range. This concludes the proof.

\subsection{Optimality of Theorem \ref{thm:bvb}}

The theorem we have presented in the previous subsection, Theorem \ref{thm:bvb}, states that Burer-Monteiro methods succeed for almost all problems satisfying Assumptions \ref{hyp:exact}, \ref{hyp:SDP}, \ref{hyp:submersion}, provided that
\begin{equation}\label{eq:constraint_p}
\frac{p(p+1)}{2} >m.
\end{equation}
The assumptions can be modified so as to include more problems \citep*{bhojanapalli_boumal_jain_netrapalli}. In addition, if one applies a slight perturbation to the cost matrix $C$ before running a suitable algorithm, it is possible to guarantee that Burer-Monteiro methods succeed at solving all (and not almost all) problems satisfying the assumptions with high probability, and with a running time polynomial in the precision \citep*{pumir_jelassi_boumal,cifuentes_moitra}.

This means that, when Inequality \eqref{eq:constraint_p} holds, Burer-Monteiro methods perform well, and are also relatively well understood from a theoretical point of view. Unfortunately, Inequality \eqref{eq:constraint_p} is quite disappointing. Indeed, it can only be satisfied if $p\geq \sqrt{2m}+o(1)$. Let us recall that Burer-Monteiro methods are applicable as soon as $p\geq \mathrm{rank}(X^s)$, and that they numerically seem to almost always succeed even if $p$ is only slightly larger than $\mathrm{rank}(X^s)$ (Paragraph \ref{sss:burer_monteiro_num}). In practice, since the computational cost increases with $p$, the most frequently used values of $p$ are therefore of order $\mathrm{rank}(X^s)$, and not of order $\sqrt{2m}$\footnote{Typically, $\mathrm{rank}(X^s)\ll \sqrt{2m}$. In many interesting applications, $\mathrm{rank}(X^s)=O(1)$ while $\sqrt{2m}=O(\sqrt{n})$.} Consequently, Inequality \eqref{eq:constraint_p} is rarely satisfied in concrete situations.

It is therefore of importance to determine whether Condition \eqref{eq:constraint_p} is optimal (that is, necessary in order for Theorem \ref{thm:bvb} to hold) or whether it can be improved. This is the question considered in \citep*{waldspurger_waters}, and the answer is that the condition is essentially optimal.
\begin{thm}\label{thm:optimality}
  Let $\mathcal{A},b$ be fixed. We assume that Assumptions \ref{hyp:exact}, \ref{hyp:SDP} and \ref{hyp:submersion} are verified, as well as a ``minimal intersection'' assumption (see below).

  Let us define $r_0=\min\{\mathrm{rank}(X),X\in\mathcal{E}_{SDP}\}$.

  For all $p$ such that
  \begin{equation}\label{eq:constraint_p_r0}
    \frac{p(p+1)}{2}+pr_0\leq m,
  \end{equation}
  there exists a set $\mathcal{C}\subset\mathcal{S}_n(\R)$ with non-zero Lebesgue measure such that, for all $C\in\mathcal{C}$,
  \begin{enumerate}
  \item \label{it:optimality_conc1} $\mathrm{rank}(X^s)=r_0$;
  \item \label{it:optimality_conc2} Problem \eqref{eq:SDP_fact} has a second-order critical point which is not a global minimizer.
  \end{enumerate}
\end{thm}

Why does Theorem \ref{thm:optimality} imply, as asserted above, that Condition \eqref{eq:constraint_p} is essentially necessary for Theorem \ref{thm:bvb} to hold? Let us remember that, in most applications, the matrices to be recovered have very small rank, typically of order $1$, hence $r_0=O(1)$. Therefore, Inequality \eqref{eq:constraint_p_r0} is verified when $p\leq \sqrt{2m}+o(1)$: up to this $o(1)$, this is the exact negation of Inequality \eqref{eq:constraint_p}. As a consequence, we can rephrase Theorem \ref{thm:optimality} as: ``When Inequality \eqref{eq:constraint_p} does not hold, Theorem \ref{thm:bvb} is false: there exists a non-negligible set of cost matrices $C$ for which Problem \eqref{eq:SDP_fact} has bad critical points, and this is even if we restrict ourselves to problem instances where the solution has the smallest possible rank $r_0$''.

We do not provide the exact definition of the ``minimal intersection'' assumption here. We invite the curious reader to look for it in \citep*{waldspurger_waters}. It consists in requiring that the intersection between two specific subspaces of $\R^{n\times p}$ is as small as possible. I do not have a satisfactory geometric interpretation of this condition, but it is necessary for the proof of Theorem \ref{thm:optimality}. Fortunately, there are good reasons to think that it is ``generically'' satisfied and it is at least satisfied in all applications discussed in \citep*{waldspurger_waters}. It is therefore not a significant restriction to the applicability of Theorem \ref{thm:optimality}.

\subsubsection{Idea of proof for Theorem \ref{thm:optimality}}

The proof is in two parts.
\begin{enumerate}
\item First part: we show that, when a cost matrix $C$ satisfies Properties \ref{it:optimality_conc1} and \ref{it:optimality_conc2} of Theorem \ref{thm:optimality} and some additional ``non-degeneracy'' conditions, then all cost matrices in some neighborhood of $C$ satisfy Properties \ref{it:optimality_conc1} and \ref{it:optimality_conc2} as well.

  As a consequence, in order to prove that there exists a set with non-zero Lebesgue measure of cost matrices with Properties \ref{it:optimality_conc1} and \ref{it:optimality_conc2}, it suffices to show that there exists \textit{one} cost matrix $C$ satisfying these two properties as well as the non-degeneracy conditions.
\item Second part: we show the existence of this one cost matrix.
\end{enumerate}

We do not describe the first part of the proof, which relies on standard arguments from differential geometry. Let us focus on the second part. We will present the main arguments of this part, and try to explain where the condition $\frac{p(p+1)}{2}+pr_0 \leq m$ comes from. For simplicity, we ignore the non-degeneracy conditions. We must therefore only explain how to construct $C$ satisfying Properties \ref{it:optimality_conc1} and \ref{it:optimality_conc2}.

Let us fix $X_0\in\mathcal{E}_{SDP}$ with rank $r_0$ and $V\in\mathcal{M}_p$. We are going to construct $C$ such that\footnote{Slightly surprinsingly, a cost matrix $C$ satisfying the desired properties exists for almost all choices of $X_0\in\mathcal{E}_{SDP},V\in\mathcal{M}_p$ such that $\mathrm{rank}(X_0)=r_0$.}
\begin{itemize}
\item the solution $X^s$ of Problem \eqref{eq:SDP} is $X_0$, hence $\mathrm{rank}(X^s)=r_0$ and Property \ref{it:optimality_conc1} holds;
\item $V$ is a second-order critical point of Problem \eqref{eq:SDP_fact}, but not a global minimizer, hence Property \ref{it:optimality_conc2} also holds.
\end{itemize}

\paragraph{First step: analytic formulations}
We first rewrite these properties under a more analytic form.

A sufficient (and almost necessary) condition for the first property to hold is given by the classical duality theory of semidefinite programs. We state it in the following proposition.
\begin{prop}\label{prop:reformulation_prop1} If there exist $g_1\in\R^m,C_1\in\mathcal{S}_n(\R)$ such that
\begin{gather*}
  C = \mathcal{A}^*(g_1) + C_1, \\
  C_1 X_0 = 0, \\
  C_1 \succeq 0,\\
  \mathrm{rank}(C_1)=n-r_0,
\end{gather*}
then $X_0$ is the unique solution of Problem \eqref{eq:SDP}.
\end{prop}

The second property can be reformulated using the explicit formulas of $\nabla_{\mathcal{M}_p}f_C$ and $\nabla_{\mathcal{M}_p}^2f_C$.
\begin{prop}\label{prop:reformulation_prop2}
$V$ is a second-order critical point of Problem \eqref{eq:SDP_fact} if and only if
there exist $g_2\in\R^m,C_2\in\mathcal{S}_n(\R)$ such that
\begin{gather*}
  C = \mathcal{A}^*(g_2) + C_2, \\
  C_2V = 0, \\
  \forall \dot{V}\in T_V\mathcal{M}_p,\quad
  \scal{C_2}{\dot V \dot{V}^T} \geq 0.  
\end{gather*}
\end{prop}

From Propositions \ref{prop:reformulation_prop1} and \ref{prop:reformulation_prop2}, in order to construct $C$ as desired, we simply need to find $g_1,g_2,C_1,C_2$ such that
\begin{gather}
  \mathcal{A}^*(g_1) + C_1 = \mathcal{A}^*(g_2) + C_2, \label{eq:C_cond1} \\
  C_1 X_0 = 0, \\
  C_1 \succeq 0, \\
  \mathrm{rank}(C_1)=n-r_0, \\
  C_2 V = 0, \label{eq:C_cond5} \\
  \forall \dot{V} \in T_V\mathcal{M}_p,\quad
  \scal{C_2}{\dot{V}\dot{V}^T} \geq 0. \label{eq:C_cond6}
\end{gather}
Indeed, if we find such $g_1,g_2,C_1,C_2$, the matrix $C=\mathcal{A}^*(g_1)+C_1$ satisfies Properties \ref{it:optimality_conc1} and \ref{it:optimality_conc2}.

\paragraph{Second step: construction of $\boldsymbol{g_1,g_2,C_1,C_2}$}
Without loss of generality, we can set $g_2=0$ and construct $g_1,C_1,C_2$ only.

One can also show that, if $g_1,C_1,C_2$ satisfy Properties \eqref{eq:C_cond1} to \eqref{eq:C_cond5}, it is possible to define $\tilde{C}_1,\tilde{C}_2$, by adding a suitable semidefinite positive matrix to $C_1$ and $C_2$, such that $g_1,\tilde{C}_1,\tilde{C}_2$ also satisfy
Properties \eqref{eq:C_cond1} to \eqref{eq:C_cond5}, and Property \eqref{eq:C_cond6} as well. Consequently, it is enough to construct $g_1,C_1,C_2$ with Properties \eqref{eq:C_cond1} to \eqref{eq:C_cond5}.

Let us explain the construction in the simplified setting where
\begin{gather*}
  \mathrm{Range}(X_0) = \{(x_1,\dots,x_{r_0},0,\dots,0)\mbox{ with }x_1,\dots,x_{r_0}\in\R\}, \\
  \mathrm{Range}(V) = \{(0,\dots,0,x_{r_0+1},\dots,x_{r_0+p},0,\dots,0)\mbox{ with }
  x_{r_0+1},\dots,x_{r_0+p}\in\R\}.
\end{gather*}

In this setting, we have the following proposition.
\begin{prop}\label{prop:condition_on_g1}
  For any $g_1\in\R^m$, the following properties are equivalent.
  \begin{enumerate}
  \item There exist $C_1,C_2\in\mathcal{S}_n(\R)$ such that $g_1,C_1,C_2$ satisfy Properties \eqref{eq:C_cond1} to \eqref{eq:C_cond5}.
  \item $\mathcal{A}^*(g_1)$ admits a block-decomposition of the form
    \begin{equation}\label{eq:Astar_g1_decomposition}
      \mathcal{A}^*(g_1) = \left(
        \begin{smallmatrix}G_1&0&G_2 \\ 0&G_3&G_4\\G_2^T&G_4^T&G_5\end{smallmatrix}
      \right),
    \end{equation}
    for some $G_1\in\mathcal{S}_{r_0}(\R)$, $G_2\in\R^{r_0\times (n-(r_0+p))}$,
    $G_3\in\mathcal{S}_p(\R)$, $G_4\in\R^{p\times (n-(r_0+p))}$, $G_5\in\mathcal{S}_{n-(r_0+p)}(\R)$, such that $G_3 \prec 0$.
  \end{enumerate}
\end{prop}

We observe that the mapping
\begin{equation*}
  \begin{array}{rccccc}
    \Phi:&\R^m&\to&\R^{r_0\times p}&\times&\mathcal{S}_p(\R) \\
         & g &\to& \Big((\mathcal{A}^*(g)_{k,l})_{\substack{1\leq k\leq r_0 \\ r_0<l\leq r_0+p}}
         &,&(\mathcal{A}^*(g)_{k,l})_{r_0<k,l\leq r_0+p}\Big)
  \end{array}
\end{equation*}
is a linear map between a vector space of dimension $m$ and a vector space of dimension
\begin{equation*}
  \frac{p(p+1)}{2} + r_0 p \leq m.
\end{equation*}
We can thus expect that it is generically surjective. And actually, using the minimal intersection property, we can rigorously guarantee that it is surjective.

Let us now fix any $H\prec 0$ in $\mathcal{S}_p(\R)$. Since $\Phi$ is surjective, there exists $g_1$ such that $\Phi(g_1)=(0,H)$. Let us fix one such $g_1$. Then the matrix $\mathcal{A}^*(g_1)$ admits a block-decomposition of the form \eqref{eq:Astar_g1_decomposition} (with $G_3=H$). From Proposition \ref{prop:condition_on_g1}, this is enough to ensure the existence of $C_1,C_2$ such that $g_1,C_1,C_2$ satisfy Properties \eqref{eq:C_cond1} to \eqref{eq:C_cond5}.

\subsection{Summary and open questions}

To summarize, we can distinguish three regimes regarding the behavior of the Burer-Monteiro heuristic.
\begin{itemize}
\item When the factorization rank $p$ is equal to, or barely larger than $\mathrm{rank}(X^s)$, it is numerically observed (see Figure \ref{fig:bm_phasecut}, for instance) that Burer-Monteiro methods succeed at solving some problems, but that there are also situations of practical interest in which they fail.
\item When $p$ is larger than $\mathrm{rank}(X^s)$ but smaller than $\sqrt{2m}+o(1)$, Burer-Monteiro methods seem to correctly find a global minimum in almost all situations of practical interest, but no theoretical explanation of this phenomenon exists.
\item When $p$ is larger than $\sqrt{2m}+o(1)$, Burer-Monteiro methods almost always converge to a global minimum and satisfactory theoretical guarantees exist (Theorem \ref{thm:bvb}).
\end{itemize}

In practice, computational efficiency commands to choose $p$ as small as possible among the values which allow to find a global minimum. The most interesting regime is therefore the second one, when $\mathrm{rank}(X^s) < p \leq \sqrt{2m} + o(1)$. Thus, a major open question is:
\begin{center}
  How can we explain the good numerical behavior of Burer-Monteiro heuristics

  in this regime?
\end{center}

More precisely, it follows from Theorem \ref{thm:optimality} that, when $p\leq \sqrt{2m}+o(1)$, there exists a non-zero Lebesgue measure set of cost matrices $C$ for which Problem \eqref{eq:SDP_fact} has bad second-order critical points. How come that we do not seem to encounter these matrices in numerical experiments? Is it because problematic cost matrices form a subset of $\mathcal{S}_n(\R)$ with an extremely small, although non-zero, volume, which makes it extremely unlikely to encounter one of them in an experiment? Is it that we encounter them, but that the attraction basin of bad critical points is very small, so that convergence to a global minimum occurs despite the existence of second-order critical points?

Other open questions have to do with more concrete algorithmic aspects of Burer-Monteiro methods. For instance, how to deal with problems for which Assumption \ref{hyp:SDP} or \ref{hyp:submersion} does not hold? For Assumption \ref{hyp:SDP}, this has already been studied in \citep*{bhojanapalli_boumal_jain_netrapalli}. Independently, which Riemannian optimization algorithms are best suited to which problems? This is a crucial issue in applications, but probably a delicate one. A notable source of difficulties is the possible ill-conditioning of Problem \eqref{eq:SDP_fact} close to the solution. About this point, the interested reader can refer to \citep*{tong_ma_chi} and references therein.

% \bibliographystyle{plainnat}
% \bibliography{../../../biblio/bib_articles.bib,../../../biblio/bib_proceedings.bib,../../../biblio/bib_livres.bib,../../../biblio/bib_misc.bib}

\appendix

\section{Proof of Proposition \ref{prop:fixed_point}\label{s:fixed_point}}

For all $k$,
  \begin{align*}
    z^{obj *}Cz^{obj}
    & = \sum_{l,l' \ne k}C_{l,l'}\overline{z^{obj}_l} z^{obj}_{l'}
      + 2 \Re\left(\overline{z^{obj}_k} \sum_{l'\ne k}C_{k,l'}z^{obj}_{l'} \right)
      + 1.
  \end{align*}
  As $z^{obj *}Cz^{obj}=\max_{|z_1|=\dots=|z_n|=1}z^* Cz$, we must have
  \begin{align*}
    \sum_{l,l' \ne k}C_{l,l'}\overline{z^{obj}_l} z^{obj}_{l'}
    & + 2 \Re\left(\overline{z^{obj}_k} \sum_{l'\ne k}C_{k,l'}z^{obj}_{l'} \right)
      + 1 \\
    & = \max_{\theta\in\R} \left(\sum_{l,l' \ne k}C_{l,l'}\overline{z^{obj}_l} z^{obj}_{l'}
      + 2 \Re\left(e^{-i\theta} \sum_{l'\ne k}C_{k,l'}z^{obj}_{l'} \right)
      + 1\right),
  \end{align*}
  which is equivalent to the existence of some $\lambda_k\in\R^+$ such that
  \begin{equation*}
    \sum_{l'\ne k}C_{k,l'}z_{l'}^{obj} = \lambda_k z_k^{obj}.
  \end{equation*}
  Therefore,
  \begin{equation*}
    (Cz^{obj})_k = z^{obj}_k + \sum_{l'\ne k}C_{k,l'}z_{l'}^{obj} = (1+\lambda_k)z^{obj}_k,
  \end{equation*}
  which implies
  \begin{align*}
    \mathcal{P}(C z^{obj})_k
    & = \frac{(1+\lambda_k)z^{obj}_k}{\left|(1+\lambda_k)z^{obj}_k\right|} \\
    & = z^{obj}_k.
  \end{align*}

\end{document}